\documentclass[12pt]{amsart}
\usepackage{amsfonts,amsmath,amssymb}

\newtheorem{Lem}{Lemma}[section]
\newtheorem{Prop}[Lem]{Proposition}
\newtheorem{Cor}[Lem]{Corollary}

\newtheorem{Theorem}[Lem]{Theorem}
\theoremstyle{definition}
\newtheorem{Def}[Lem]{Definition}
\newtheorem{Rem}[Lem]{Remark}

\newtheorem{Ex}[Lem]{Example}

\newcommand{\pf}{\medskip\noindent{\sc Proof: }}
\newcommand{\epf}{$\Box$}

\newcommand\Alg{\operatorname{Alg}}

\newcommand\id{{\operatorname{id}}}

\newcommand\cop{{\operatorname{cop}}}
\newcommand\ad{{\operatorname{ad}}}

\newcommand\ord{{\operatorname{ord}}}

\newcommand\Irr{{\operatorname{Irr}}}

\renewcommand\o{\otimes}

\newcommand\YDG{^{\Gamma}_{\Gamma}\mathcal{YD}}
\newcommand\YDL{^{\Lambda}_{\Lambda}\mathcal{YD}}

\newcommand\YDI{^{\mathbb{Z}[I]}_{\mathbb{Z}[I]}\mathcal{YD}}

\newcommand\D{\mathcal{D}}

\newcommand\G{\Gamma}

\newcommand\sw[1]{{}_{(#1)}}

\title[Irreducible representations]{On the irreducible representations of a class of pointed Hopf algebras}
\author{David E. Radford}
\address{University of Illinois at Chicago \\
Department of Mathematics, Statistics and \\ Computer Science (m/c 240) \\
801 South Morgan Street \\
Chicago, IL   60607-7045} \email{radford@uic.edu}
\author{Hans J\"{u}rgen Schneider}
\address{Mathematisches Institut \\
Ludwig-Maximilians-Universit\"{a}t M\"{u}nchen \\
Theresienstr. 39 \\
D-80333 M\"{u}nchen, Germany}
\email{Hans-Juergen.Schneider@mathematik.uni-muenchen.de}
\numberwithin{equation}{section}
\begin{document}

\date{}
\begin{abstract}
{ \small  \rm We parameterize the finite-dimensional irreducible representations of a class of pointed Hopf algebras over an algebraically closed field of characteristic zero by dominant characters. The Hopf algebras we are considering arise in the work of N. Andruskiewitsch and the second author. Special cases are the multiparameter deformations of the enveloping algebras of semisimple Lie algebras where the deforming parameters are not roots of unity and some of their finite-dimensional versions in the root of unity case.}
\end{abstract}

\maketitle

\setcounter{section}{0}
\section*{Introduction}\label{SecIntro}
In this paper we parameterize the finite-dimensional irreducible
representations of pointed Hopf algebras over an
algebraically closed field $k$ of characteristic zero
which arise in recent classification work \cite{AS2,AS3}.
These Hopf algebras are quotients of two-cocycle twists
$$H=(U \otimes A)^{\sigma},$$
where $U$ and $A$ have the form $B \# k[G]$ \phantom{a}
where $k[G]$ is the group algebra of an abelian
group $G$ over $k$ and $B$ is a left $k[G]$-module algebra
and a left $k[G]$-comodule coalgebra \cite{RP}. More
precisely $B= \mathfrak{B}(X)$ is the Nichols algebra of a
finite-dimensional Yetter-Drinfeld module over the group algebra
$k[G]$. See \cite[Section 2]{ASSurvey} for a discussion of Nichols
algebras and Yetter-Drinfeld modules in general.

Thus
$$U=\mathfrak{B}(W){\#}k[\Lambda] \mbox{$\;\;$ and $\;\;$} A = \mathfrak{B}(V){\#}k[\G],$$
where $\Lambda$ and $\G$ are abelian groups, and $W$ and $V$ are
Yetter-Drinfeld modules over $\Lambda$ and $\G$
respectively. We assume that $W$ and $V$ are direct sums of
one-dimensional Yetter-Drinfeld modules. The 2-cocycle $\sigma :
(U \otimes A) \otimes (U \otimes A) \to k$ is given in terms of a
bilinear form $\beta : W \otimes V \to k$.
\medskip

In the Sections \ref{Sectiongeneral} and \ref{computations}
we study the set $\Irr(H)$ of isomorphism classes of
finite-dimensional irreducible $H$-modules. We assume that
the finite-dimensional irreducible $U$- and $A$-modules are
one-dimensional; by Theorem \ref{dim1} this assumption is
satisfied in the infinite-dimensional generic case, that is, when
the diagonal elements of the braiding matrix of $W$ are not roots
of unity, and when $U$ and $A$ are finite-dimensional.  Then we
have shown in \cite{Pairs} that the finite-dimensional irreducible
$H$-modules have the form $L(\rho,\chi)$, where $\rho \in
\widehat{\Lambda}$ and $\chi \in \widehat{\G}$ are certain
characters. We conclude from results of \cite{DERHJS2} and
\cite{Pairs} that there are Yetter-Drinfeld submodules $W' \subset
W$ and $V' \subset V$, and a Hopf algebra projection $H  \to H'$
defining a bijection
\begin{equation}\label{intronondeg}
\Irr(H') \xrightarrow{\cong} \Irr(H),
\end{equation}
where $H'=(\mathfrak{B}(W') \# k[\Lambda] \otimes \mathfrak{B}(V')
\# k[\G])^{\sigma'}$, and $ \sigma'$ is the restriction of
$\sigma$ such that the restricted bilinear form $\beta' : W'
\otimes V' \to k$ is {\em non-degenerate}.

In Section \ref{computations} we assume that the bilinear
form $\beta$ is non-degenerate. The main result of this section is
Theorem \ref{maingeneric} where we show that
\begin{equation}
\{(\rho,\chi) \in \widehat{\Lambda} \times \widehat{\G} \mid
(\rho,\chi) \text{ dominant}\} \to \Irr((U \otimes A)^{\sigma}),
\end{equation}
given by mapping $(\rho,\chi)$ onto the isomorphism class of
$L(\rho,\chi)$, is bijective. The assumption in Theorem
\ref{maingeneric} is that the braiding matrix  of $W$ is generic,
that is all its diagonal entries are not roots of unity, and of
finite Cartan type \cite{AS1}. The notion of a {\em dominant pair
of characters} $(\rho,\chi)$ is defined in \eqref{dominant}.
\medskip

In Section \ref{SectionCartan} we apply the
results of the first two sections to the class of pointed
Hopf algebras studied in \cite{AS3}, in particular to
the Hopf algebras of the type $U(\D,\lambda)$ or $u(\D,\lambda)$
described below. Let $(g_i)_{1 \leq i \leq \theta}$ be elements
in $\G$, let $(\chi_i)_{1 \leq i \leq \theta}$ be
elements of $\widehat{\G}$, the $k$-valued
characters of $\G$, and let $(a_{ij})_{1 \leq i,j
\leq \theta}$ \phantom{a}  be a Cartan matrix of
finite type. Then the collection
$$\mathcal{D} = \mathcal{D}(\Gamma, (g_i), (\chi_i), (a_{ij}))$$
is called a datum of finite Cartan type \cite{AS3} if
\begin{equation*}
q_{ij} q_{ji} = q_{ii}^{a_{ij}},\;q_{ii}\neq 1, \text{ with }
q_{ij} = \chi_j(g_i) \text{ for all } 1 \leq i,j \leq \theta.
\end{equation*}
$\D$ is called generic if no $q_{ii}$ is a root of unity.

Let $ \lambda = (\lambda_{ij})_{1 \le i, j \le \theta, \,
i\not\sim j}$ be a family of linking parameters for $\D$, that is,
$\lambda_{ij} \in k$ for all $1 \le i, j \le \theta, \, i\not\sim
j$, $\lambda_{ji} = - q_{ji} \lambda_{ij}$, and if $g_ig_j =1
\text{ or }\chi_i \chi_j \neq \varepsilon$,  then $\lambda_{ij}
=0$.

The Hopf algebra $U(\D,\lambda)$ is generated as an algebra by the
group $\G,$ that is, by generators of $\G$ satisfying the
relations of the group, and $x_1,\dots,x_\theta,$ with the
relations:
\begin{align*}
&\text{({\em Action of the group}) }& &gx_i g^{-1} = \chi_i(g)
x_i,  \text{ for all }
i, \text{ and all } g \in \G,&\\
&\text{({\em Serre relations}) }& &\ad_c(x_i)^{1 - a_{ij}}(x_j) = 0, \text{ for all }  i \neq j, i \sim j,&\\
&\text{({\em Linking relations}) }&& \ad_c(x_i)(x_j) = \lambda_{ij}(1 - g_ig_j), \text{ for all }i \nsim j.&\\
\intertext{The coalgebra structure is given by} &\Delta (x_i) =
g_i \o x_i + x_i \o 1, &&\Delta(g) = g \o g, \text{ for all }1
\leq i \leq \theta,  g \in \G.&
\end{align*}
In \cite{AS3} a glueing process was used to build  $U(\D,\lambda)$
inductively by adding one connected component at a time and
modding out some central group-like elements in each step. The
methods of the present paper to parameterize the
finite-dimensional irreducible modules do not apply to this
description of $U(\D,\lambda)$, since in Theorem \ref{maingeneric}
we have to assume that the finite-dimensional irreducible modules
over $U$ and $A$ are one-dimensional. In Lemma \ref{linkinggraph}
we show that the linking graph of $(\D,\lambda)$ in the generic
case has no odd cycles. Hence the linking graph is bipartite, and
we can give another description of $U(\D,\lambda)$ by one glueing
only. As a consequence Theorem \ref{maingeneric} applies and in
Theorem \ref{mainD} we obtain a bijection
\begin{equation}\label{mainintro}
\{\chi \in \widehat{\G} \mid \chi \text{ dominant }\} \to
\Irr(U(D,\lambda)),
\end{equation}
where dominant characters $\chi$ are defined in
\eqref{dominantchi}. In the proof we use the reduction to the
non-degenerate case by \eqref{intronondeg}.

The Hopf algebras $U(\D,\lambda)$ with generic braiding matrix
form a very general class of pointed Hopf algebras of finite
Gelfand-Kirillov dimension and with abelian group of group-like
elements.   In Section \ref{sectionreduced} we introduce reduced
data $\D_{red}$ of finite Cartan type and define Hopf algebras
$U(\D_{red},l)$, where $l=(l_i)_{1 \leq i \leq n}$ is a family of
non-zero scalars. The Hopf algebras $U(\D_{red},l)$ are a
reformulation of $U(\D,\lambda)$ in the non-degenerate case, that
is, when any vertex of $\{1,\dots,\theta\}$ is linked to some
other vertex. Thus to study the finite-dimensional irreducible
$U(\D,\lambda)$-modules it suffices to consider modules over
$U(\D_{red},l)$.

In Lemma \ref{EF} we describe $U(\D_{red},l)$ by the usual
generators $E_i$ and $F_i, 1 \leq i \leq n$, and  deformed Serre
relations. In the end of Section \ref{sectionreduced} we note that
the one-parameter deformation $U_q(\mathfrak{g})$ of
$U(\mathfrak{g})$, where $\mathfrak{g}$ is a semisimple Lie
algebra (see \cite{J}), Lusztig's version of the one-parameter
deformation with more general group-like elements in \cite{L}, and
the two-parameter deformations of $U(\mathfrak{gl_n})$ and
$U(\mathfrak{sl_n})$ discussed in \cite{BW} are all special cases
of $U(\D_{red},l)$. In all these cases the parametrization of the
finite-dimensional irreducible representations (of ``type 1'') by
dominant weights is a special case of \eqref{mainintro}.
\medskip

Our results on parameterization of the irreducible represenations
can be partially extended to the finite-dimensional versions of
$U(\D,\lambda)$ in \cite{Class}, where $\G$ is a finite abelian
group. We assume that the root vector relations in \cite{Class}
are all $0$. The finite-dimensional Hopf algebras $u(\D,\lambda)$
defined in Section \ref{finite} generalize Lusztig's Frobenius
kernels. The linking graph of $u(\D,\lambda)$ is bipartite if the
Cartan matrix is simply laced but not in general. The
finite-dimensional version of Theorem \ref{maingeneric} applies to
$u(\D,\lambda)$ with bipartite linking graph, and in Theorem
\ref{finitesimple} we obtain a bijection
\begin{equation}
\widehat{\G} \to \Irr(u(\D,\lambda)).
\end{equation}

Throughout this paper $k$ is an algebraically closed field
of characteristic zero.

\section{Preliminaries and general theorems}\label{Sectiongeneral}

We recall and reformulate some general results from \cite{DERHJS2}
and \cite{Pairs} for the class of Hopf algebras we are
considering. We refer to \cite[Section 2]{ASSurvey} for a
discussion of braided Hopf algebras, Nichols algebras and
Yetter-Drinfeld modules in general, and to \cite{SweedlerBook} and
\cite{Mont} for Hopf algebra theory.

Let $\Gamma$ be an abelian group. A Yetter-Drinfeld module over
$k[\G]$ can be described as a $\Gamma$-graded vector space which
is a $\Gamma$-module such that all  $g$-homogeneous components,
where $g \in \Gamma$, are stable under the $\Gamma$-action. We
denote the category of Yetter-Drinfeld modules over $k[\G]$ by
${}_{k[\G]}^{k[\G]}{\mathcal YD} = {\YDG}.$

For $X\in {\YDG},$ $g \in \Gamma,$ and $\chi \in \widehat{\Gamma}$
we define $$X_g = \{x \in X \mid \delta(x) = g \otimes x\}$$ and
$$X_g^{\chi} = \{ x \in X_g \mid h \cdot x = \chi(h)x \text{ for all } h \in \Gamma \}.$$
The category $\YDG$ is braided and for $X,Y \in {\YDG}$
$$c: X \otimes Y \to Y \otimes X, x \otimes y \mapsto g \cdot y \otimes x, x \in X_g,y \in Y,$$
defines the braiding on $X \otimes Y$.

Let $ \theta \geq 1, g_1,\dots,g_{\theta} \in \G,\chi_1,\dots,
\chi_{\theta} \in \widehat{\G}$, and 
$$q_{ij} = \chi_j(g_i), 1 \leq
i,j \leq \theta.$$ 
Let $X \in {\YDG}$ be the object with basis $x_i \in
X_{g_i}^{\chi_i},1 \leq i \leq \theta$. The braiding on the vector space $X$ is given by 
$$c(x_i \otimes x_j)= q_{ij}x_j \otimes x_i,1
\leq i,j \leq \theta.$$

A braided graded Hopf algebra
$$R = \oplus_{n \geq 0} R(n) \text{ in }\YDG$$
is a {\em Nichols algebra} of $X$ if $X \cong R(1)$ in $\YDG$, $R$
is connected, that is, $k \cong R(0)$,  and if
\begin{enumerate}
\item $R(1)$ consists of all the primitive elements of $R$ and 
\item as an algebra $R$ is generated by $R(1)$.
\end{enumerate}
The Nichols algebra of $X$ exists and is unique up to isomorphism. We denote the Nichols algebra of $X$
by $ \mathfrak{B}(X)$. As an algebra and coalgebra,
$\mathfrak{B}(X)$ only depends on the braided vector space
$(X,c)$. The structure of $\mathfrak{B}(X)$ is known if the
braiding of $X$ is related to a Cartan matrix of finite type in
the following way.

Let $(a_{ij})_{1 \leq i,j \leq \theta})\phantom{a}$ be a Cartan matrix of
finite type. Then the collection
$$\mathcal{D} = \mathcal{D}(\Gamma, (g_i)_{1 \leq i \leq \theta}, (\chi_i)_{1 \leq i \leq \theta}, (a_{ij})_{1 \leq i,j \leq \theta})$$
is called a {\em datum of finite Cartan type} and the {\em
braiding of $X$ is  of finite Cartan type} \cite{AS3} if
\begin{equation}\label{Cartan}
q_{ij} q_{ji} = q_{ii}^{a_{ij}},\;q_{ii}\neq 1, \text{ where }
q_{ij} = \chi_j(g_i) \text{ for all } 1 \leq i,j \leq \theta.
\end{equation}
The smash product $\mathfrak{B}(X) \# k[\G]$ is a  Hopf
algebra over the field $k$ where
\begin{align}
gx_ig^{-1} &= \chi_i(g)x_i ,\text{ and }\label{actionU}\\
\Delta(x_i)&=g_i \otimes x_i + x_i \otimes 1,\; \Delta(g) = g
\otimes g\label{comultU}
\end{align}
for all $g \in \G, 1 \leq i \leq \theta$. Here and in the
following we
 identify $x \in \mathfrak{B}(X)$ with $x\# 1$, and for $h \in k[\G]$ we identify $h$ with $1 \#h$ in the smash product. Thus we write
\begin{align*}
xh= (x\#1)(1\#h) = x \#h,\text{ and } hx=(1\#h)(1\#b) =
(h\sw1\cdot b)h\sw2.
\end{align*}

Suppose that the braiding of $X$ is of finite Cartan type with
Cartan matrix $(a_{ij})$. For all $ 1 \leq i,j \leq n$ we write $i
\sim j$ if $i$ and $j$ are in the same connected component of the
Dynkin diagram of $(a_{ij})$. We remark that $\sim$ is an equivalence relation and the connected components of $I = \{1,2,\dots,\theta\}$ are its equivalence classes, the set of which we denote by $\mathcal{X}$.

By \cite[Lemma 2.3]{Class} there are $d_i \in \{1,2,3\},1 \leq i
\leq \theta$, and $q_{J} \in k$ for all $J \in \mathcal{X}$, such that 
\begin{equation}\label{DJ}
q_{ii}=q_J^{2d_i}, d_ia_{ij} =d_ja_{ji} \text{ for all }J \in
\mathcal{X},i,j \in J.
\end{equation}
We define
\begin{equation}\label{q'}
q'_{ij} =\begin{cases} (q_J)^{d_ia_{ij}},& \text{ if } i,j \in J,J \in \mathcal{X},\\
0,& \text{ if } i \nsim j.
\end{cases}
\end{equation}
Then
\begin{equation}\label{twist}
q_{ij}q_{ji} = q'_{ij}q'_{ji},q_{ii} = q'_{ii} \text{ for all }1
\leq i,j \leq \theta. \end{equation} Thus in the terminology of
\cite{AS3}, the braiding of $X$ is twist equivalent to the
braiding of Drinfeld-Jimbo type given by $(q'_{ij})$.

The datum $\D$ and the braiding of $X$ are called {\em generic}
\cite{AS3} if
\begin{equation}\label{generic}
q_{ii} \text{ is not a root of  unity  for all }1 \leq i \leq n.
\end{equation}
If $\D$ is a generic Cartan datum of finite type, then it is shown
in \cite{AS3} using results of Lusztig \cite{L}  and Rosso
\cite{R}  that
\begin{equation*}
\mathfrak{B}(X)\cong k\langle x_1,\dots,x_{\theta} \mid
(\ad_c(x_i))^{1-a_{ij}}(x_j)=0, 1 \leq i,j \leq \theta, i\neq j
\rangle,
\end{equation*}
where $\ad_c$ denotes the braided adjoint action. Hence if $y =
x_{j_1} \cdots x_{j_s}$, and $1 \leq i,j_1,\dots,j_s \leq \theta,s
\geq 1$, then
$$(\ad_c(x_i))(y) = xy - q_{ij_1} \cdots q_{ij_s}yx_i.$$
In particular, if $\mathfrak{g}$ is the semisimple Lie algebra
with Cartan matrix $(a_{ij})$, $q \in k$ is not a root of unity,
and the braiding is given by \eqref{q'}, where $q_J=q$ for all
connected components $J$, then $\mathfrak{B}(X) \cong
U_q^-(\mathfrak{g})$.

\begin{Theorem}\label{dim1}
Let
$$\mathcal{D} = \mathcal{D}(\Gamma, (g_i)_{1 \leq i \leq \theta}, (\chi_i)_{1 \leq i \leq \theta}, (a_{ij})_{1 \leq i,j \leq \theta})$$
be a datum of finite Cartan type, and let $X \in {\YDG}$ be the object with basis
$x_i \in X_{g_i}^{\chi_i},1 \leq i \leq \theta$. Then all
finite-dimensional simple $\mathfrak{B}(X) \#k[\G]$-modules are
one-dimensional if
\begin{enumerate}
\item $\D$ is generic, or 
\item The group $\G$ is finite, for
all $1 \leq i \leq \theta,$ $ \ord(q_{ii})$ is odd, and
$\ord(q_{ii})$  is prime to $3$ if $i$ is in a component $G_2$ of the
Dynkin diagram of $(a_{ij})$.
\end{enumerate}
\end{Theorem}
\pf In Case (1) the theorem is a special case of \cite[Theorem
5.6]{Pairs}. In Case (2), $\mathfrak{B}(X)$ is finite-dimensional
by \cite[Theorem 1.1]{AS1}, \cite[Theorem 5.1]{Class}. Since $\mathfrak{B}(X)$ is
$\mathbb{N}$-graded it follows that the augmentation ideal
$\mathfrak{B}(X)^+$ is nilpotent. Let $U= \mathfrak{B}(X) \#
k[\G]$. Since 
$$U\mathfrak{B}(X)^+=\mathfrak{B}(X)^+U, \text{ and }
U/U\mathfrak{B}(X)^+ \cong k[\G],$$ 
the simple $U$-modules are
simple $k[\G]$-modules, hence one-dimensional. \epf

\medskip

For the rest of this section we fix abelian groups $\Lambda$ and
$\Gamma$, integers $n, m \geq 1$, elements $z_1, \ldots, z_n \in
\Lambda$ and $g_1,\ldots, g_m \in \Gamma,$ and nontrivial
characters $\eta_1, \ldots, \eta_n \in \widehat{\Lambda}$ and
$\chi_1,\ldots, \chi_m \in \widehat{\Gamma}.$

Let $W \in {\YDL}$ have basis $u_i \in W_{z_i}^{\eta_i}, 1 \leq i
\leq n$ and let $V \in {\YDG}$ have basis $a_j \in V_{g_j}^{\chi_j},
1 \leq j \leq m.$ Let $U=\mathfrak{B}(W){\#}k[\Lambda]$ and $A =
\mathfrak{B}(V){\#}k[\G].$ Note that the restriction maps
\begin{equation}\label{injectiverestriction}
\Alg(\mathfrak{B}(W) \# k[\Lambda],k) \to \widehat{\Lambda},\;
\Alg(\mathfrak{B}(V) \# k[\G],k) \to  \widehat{\Gamma},
\end{equation}
are bijective since the characters $\eta_i,\chi_j$ are all
non-trivial.

We denote the preimage of a character $\psi \in \widehat{\Lambda}$
or $\widehat{\G}$ by $\widetilde{\psi}$. Thus if $\psi \in
\widehat{\G}$, then $\widetilde{\psi}(a_j)=0$, and
$\widetilde{\psi}(g)=\psi(g)$ for all $1 \leq j \leq m,g \in \G$.

The following results hold in this general context; we do not
assume that the braidings are of finite Cartan type.

The next theorem describes how to define Hopf algebras maps from
$U$ to $A^{o\,\cop}$, the $\cop$-version of the Hopf dual of $A$.
The proof is based on the universal properties of the Nichols
algebra.

\begin{Theorem}\label{groupcase}
In addition to the above let $\varphi : \Lambda \to
\widehat{\Gamma}$ be a group homomorphism, $s : \{1,\dots,n\} \to
\{1,\dots,m\}$ be  a function, and let $l_1, \ldots, l_n \in k$.
Let $I' = \{1 \leq i \leq n \mid l_i \neq 0\}$. Assume that for
all $i \in I'$ and $z \in \Lambda$
\begin{equation}
\varphi(z_i) = \chi^{-1}_{s(i)}\; \mbox{and} \; \eta_i(z) =
\varphi(z)(g_{s(i)}).  \label{EqVarPhI2}
\end{equation}
Then  there is a Hopf algebra map $\Phi : U \to A^{o\,\cop}$ such
that
\begin{align*}
\Phi(z) &= \widetilde{\varphi(z)} \text{ for all }z \in \Lambda,\\
\Phi(u_i)(g)&=0 \text{ for all }1 \leq i \leq n,g \in \G,\\
\Phi(u_i)(a_j) &= \delta_{s(i),j} l_i \text{ for all } 1 \leq i
\leq n,1 \leq j \leq m.
\end{align*}
\end{Theorem}
\pf This is shown in \cite[Corollary 9.1]{DERHJS2}. \epf

By \cite{DoiTak} any Hopf algebra map $\Phi : U \to A^{o\,\cop}$
determines a convolution invertible map $\tau : U \otimes A \to k$
by
\begin{align}
\tau(u,a)&=\Phi(u)(a),\label{deftau}\\
\intertext{and a 2-cocycle $\sigma : (U \otimes A) \otimes (U
\otimes A) \to k$ by} \sigma(u \otimes a,u' \otimes a')&=
\varepsilon(a) \tau(u',a) \varepsilon(a'),\label{defsigma}
\end{align}
for all $u,u' \in U,a,a' \in A$. Let
$$H=(U \otimes A)^{\sigma}.$$
Recall that $(U \otimes A)^{\sigma}$ is a Hopf algebra with
coincides with $U \otimes A$ as a coalgebra with componentwise
comultiplication and whose algebra structure is defined by
\begin{align}\label{Eq2CocyleMult}
(u{\otimes}a)(u'{\otimes}a') &= \sigma(h\sw1,h'\sw1)h\sw2h'\sw2 \sigma^{-1}(h\sw3,h'\sw3)\\
&=u\tau(u'_{(1)},
a_{(1)})u'_{(2)}{\otimes}a_{(2)}\tau^{-1}(u'_{(3)},
a_{(3)})a'\notag
\end{align}
for all $u, u' \in U,a, a' \in A$, and $h=u \otimes a,h'=u'\otimes
a'$. Note that $\tau^{-1}(u,a) = \tau(S(u),a)$ for all $u \in U,
a\in A$. We view $U$ and $A$ as subalgebras of $H$ by the
embeddings $u \mapsto u \otimes 1$ and $a \mapsto 1 \otimes a$.
Since $(u \o 1)(1 \o a)=u \o a$ we have $ua = u \o a$ under our identifications of $U$ and $A$ as subalgebras of $U \o A$.

Our main concern are the irreducible finite-dimensional
representations of $H$. Let $\rho \in \widehat{\Lambda}, \chi \in
\widehat{\G}$ be characters. We define  the left $H$-module
$U_{\chi}$  with underlying vector space $U$ by
\begin{align}\label{actionchi}
(u'a)\cdot_{\chi} u &= u' \tau(u\sw1,a\sw1)u\sw2
\widetilde{\chi}(a\sw2)\tau^{-1}(u\sw3,a\sw3)
\end{align}
for all $u,u'\in U, a \in A$. Thus $U_{\chi}\cong H \otimes _A
k_{\chi},$ where $k_{\chi}$ is the one-dimensional $A$-module
defined by $\widetilde{\chi}$.

Let $I(\rho,\chi)$ be the largest $H$-submodule of $U_{\chi}$
which is contained in the kernel of $\widetilde{\rho}: U =
\mathfrak{B}(W) \# k[\Lambda] \to k$. We define
$$L(\rho,\chi) = U_{\chi}/I(\rho,\chi).$$
Then $L(\rho,\chi)$ is a cyclic left $H$-module with generator
$$m_{\rho,\chi}=\text{ residue class of }1.$$

We also need another description of $L(\rho,\chi)$ as a subspace
of the dual $A^*$.  We define  the right $H$-module $A_{\rho}$
with underlying vector space $A$ by
\begin{equation}\label{actionrho11}
a \cdot_{\rho}(ua')= \tau(u\sw1,a\sw1)\rho(u\sw2)
a\sw2a'\tau^{-1}(u\sw3,a\sw3)
\end{equation}
for all $u\in U, a,a'\in A$ . Then the dual $(A_{\rho})^*$ is a
left $H$-module by
\begin{equation}\label{actionrho}
((ua')\cdot_{\rho} p)(a)= p(a \cdot_{\rho}(ua'))
\end{equation}
for all $u\in U, a,a'\in A$ and $p \in A^*$.

Note that by \eqref{actionrho}
\begin{equation}
u\cdot_{\rho}p = \Phi(u\sw1) \rho(u\sw2)p\Phi(S(u\sw3)) \text{ for
all } u \in U,p \in A^*.\label{actionrho1}
\end{equation}

By \cite[part c) of Proposition 2.6]{Pairs}  there is a left
$H$-isomorphism
\begin{equation}\label{presentationL}
U_{\chi}/I(\rho,\chi) \cong U \cdot_{\rho} \widetilde{\chi},
\text{ with } m_{\rho,\chi} \mapsto \widetilde{\chi}.
\end{equation}
For any algebra $R$ we denote by $\Irr(R)$ the set of isomorphism classes of all finite-dimensional left $R$-modules.
\begin{Theorem}\label{modules}
In the situation of Theorem \ref{groupcase} assume that all
finite-dimensional simple $U$- and $A$-modules are
one-dimensional. Let $\sigma$ be the 2-cocycle determined by
$\Phi$ and $H=(U \otimes A)^{\sigma}$.
\begin{enumerate}
\item The map $$L_{H} : \{(\rho,\chi) \mid \rho \in
\widehat{\Lambda}, \chi \in \widehat{\Gamma}, \dim(L(\rho,\chi)) <
\infty\} \to \Irr(H),$$ given by $(\rho, \chi) \mapsto [L(\rho,
\chi)]$, is a bijection. \label{bijection} \item Let $\rho \in
\widehat{\Lambda}, \chi \in \widehat{\G}$ and let $L$ be a
finite-dimensional simple left $H$-module. Then $L \cong
L(\rho,\chi)$ if and only if there is a non-zero element $m \in L$
such that $a \cdot m = \widetilde{\chi}(a) m, z \cdot m =
\rho(z)m$ for all $a \in A, z \in \Lambda$.\label{Labstract}
\end{enumerate}
\end{Theorem}
\pf By \eqref{injectiverestriction} part \eqref{bijection} follows
from Theorem \ref{groupcase} and \cite[Theorem 4.1 ]{Pairs}.

To prove \eqref{Labstract}, let $m =m_{\rho,\chi}\in
L(\rho,\chi)$. It follows from the definition of $U_{\chi}$ that
$a \cdot 1 =\widetilde{\chi}(a)1$, hence
\begin{equation}\label{condition1}
a \cdot m = \widetilde{\chi}(a)m \text{ for all }a \in A.
\end{equation}
By \eqref{actionrho1} we have $z \cdot_{\rho} \widetilde{\chi} =
\Phi(z)\rho(z)\widetilde{\chi}\Phi(z^{-1})\rho(z) =\rho(z)
\widetilde{\chi}$ for all $z \in \Lambda$, since $\Alg(A,k)$ is
commutative by \eqref{injectiverestriction}. Hence by
\eqref{presentationL}
\begin{equation}\label{condition2}
z \cdot m = \rho(z)m \text{ for all } z \in \Lambda.
\end{equation}
We have shown in \cite[Corollary 3.5,a)]{Pairs} that
$L(\rho,\chi)$ has a unique one-dimensional $A$-submodule. Thus up
to a non-zero scalar the element $m \in L(\rho,\chi)$ is uniquely
determined by \eqref{condition1}, and the claim follows from
\eqref{condition1}. \epf

Let $\varphi : R \to S$ be any ring homomorphism. Then the
restriction functor $\varphi^* : {_{S}\mathcal M} \to
{_{R}\mathcal M}$ maps a left $S$-module $M$ to itself thought of
as an $R$-module via pullback along $\varphi$. For any algebra $R$
we denote by $\Irr(R)$ the set of isomorphism classes of all
finite-dimensional left $R$-modules.

\begin{Theorem}\label{nondegenerate}
In the situation of Theorem \ref{groupcase} assume that all
finite-dimensional simple $U$- and $A$-modules are one-dimensional
and that the restriction of $s$ to $I'$ is injective.  Let
$V'\subset V$ and $W' \subset W$ be the Yetter--Drinfeld
submodules with bases $a_{s(i)}, i \in I',$ and $u_i, i \in I'.$
Let $U' = \mathfrak{B}(W') \# k \Lambda, A' = \mathfrak{B}(V') \#
k \Gamma$ and $H'=(U' \otimes A')^{\sigma'}$, where $\sigma'$ is
the restriction of $\sigma$ to $U' \otimes A'$. Then the
projections $\pi_W : W \to W', \;\pi_V : V \to V'$ define a
surjective bialgebra map $ F: H \to H'$ determined by $F|W =
\pi_W$, $F|V=\pi_V$, $F|\G = \id_{\G}$, and $F|\Lambda =
\id_{\Lambda}$. The restriction functor $F^*$ defines a bijection
\begin{equation}
F^* : \Irr(H') \to \Irr(H),
\end{equation}
and the diagram
\begin{center}
\begin{picture}(100,80)(0,-10)
\put(104,50){\vector(0,-1){43}} \put(30,60){\vector(1,0){55}}
\put(24,49.5){\vector(4,-3){60}} \put(-5,57){$\widehat{\Lambda}
\times \widehat{\G}$} \put(90,57){$\Irr(H')$}
\put(90,-5){$\Irr(H)$} \put(50,65){$L_{H'}$} \put(106,30){$F^*$}
\put(50,35){$L_H$}
\end{picture}
\end{center}
commutes.
\end{Theorem}
\pf The first part follows by \cite[Corollary 9.2]{DERHJS2}. The
diagram commutes by \cite[Proposition 2.7]{Pairs} with
$U'=\overline{U}$, $A' = \overline{A}$, $f=F|U$, and $g=F|A$.
Hence $F^*$ is bijective by Theorem \ref{modules} for $H$ and
$H'$. \epf

\medskip

\begin{Rem}\label{sinjective}
The restriction $s|I'$ of the previous theorem is injective when
the braiding matrix $(q_{ij}=\eta_j(z_i))$ of $W$ satisfies the
following condition
$$q_{ij}q_{ji} = q_{ii}^{a_{ij}} \text{ for all } i,j,$$
where $(a_{ij})$ is a Cartan matrix of finite type and for all $i$
the order of $q_{ii}$ is greater than 3. \epf\end{Rem} 
The remark
is easily justified. By \eqref{EqVarPhI2} we have 
$$q_{ij} =
\eta_j(z_i) = \varphi(z_i)(g_{s(j)}) = \chi_{s(i)}^{-1}(g_{s(j)})$$
for all $i \in I'$ and $ 1 \leq j \leq \theta$.

Assume $s(j)=s(l)$, where $j,l \in I'$. Then $q_{ij}=q_{il}$ and
$q_{ji} =q_{li}$ for all $i \in I'$. Thus $q_{ii}^{a_{ij}} =
q_{ij}q_{ji} = q_{il}q_{li} = q_{ii}^{q_{il}}.$ Since $|a_{ij} -
a_{il}| \leq 3$ and the order of $q_{ii}$ is larger than 3 we have
$a_{ij}=a_{il}$ for all $i \in I'$. Since $(a_{rs})_{r,s \in I'}$
is invertible necessarily $j=l$.

\medskip

The results of Section 9 of \cite{DERHJS2} are given in terms of
the bilinear form $\beta : W \otimes V \to k$ defined by
$\beta(u_i \otimes a_j) = l_i \delta_{s(i)j}$ for all $1 \leq i
\leq n$ and $1 \leq j \leq m$. In the preceding theorem, this form
for $W'\otimes V'$ is $\beta'=\beta|(W'\otimes V')$ and is
non-degenerate. Most applications for us will be in the context of
Theorem \ref{groupcase} and the map $s$ will be injective by
virtue of Remark \ref{sinjective}. In light of Theorem
\ref{nondegenerate} to study the $L(\rho,\chi)'s$ we may assume
that $\beta$ is non-degenerate.

\section{The Modules $L(\rho,\chi)$ in the Non-degenerate Case}\label{computations}

We continue with the notation following Theorem \ref{dim1} and in the context of Theorem \ref{groupcase}, imposing further restrictions on the hypothesis. We
assume that $\beta$ is non-degenerate, in particular that the map
$s$ is bijective. Without restriction we assume that $s=\id$.

Let $\varphi : \Lambda \to \widehat{\G}$ be a group homomorphism,
and $l_1,\dots,l_n \neq 0$. By rescaling the generators $a_i$ we
could further assume that $l_i=1$ for all $1 \leq i \leq n$;
however this might not be convenient for computations.

As a result od our assumptions
\begin{equation}\label{conditionvarphi}
\varphi (z_i) = \chi^{-1}_{i} \text{ and }\eta_i(z) =
\varphi(z)(g_i) \text{ for all } 1 \leq i \leq n, z \in \Lambda,
\end{equation}
and the Hopf algebra map $\Phi : U \to A^{o\,\cop}$ is given by
\begin{equation}
\Phi(z) = \widetilde{\varphi(z)}, \; \Phi(u_i) = \delta_i \text{
for all } z \in \Lambda, 1 \leq i \leq n,
\end{equation}
where $\delta_i : A \to k$ is the
$(\varepsilon,\widetilde{\chi_i^{-1}})$-derivation determined by
\begin{equation}\label{definitiondelta}
\delta_i(g) = 0,\;\delta_i(a_j) = \delta_{ij} l_i \text{ for all }
g \in \Gamma,1 \leq i,j \leq n.
\end{equation}

\begin{Rem}\label{Phiinjective}
The Hopf algebra map $\Phi : U \to A^{o\, cop}$ is injective since
it is injective on the primitive elements $P(\mathfrak{B}(W))=W$
by our assumption that $l_i \neq 0$ for all $i$ (see \cite[Lemma
11.01]{SweedlerBook}). \epf\end{Rem}

Let $H=(U \otimes A)^{\sigma}$, where $\sigma$ is the 2-cocycle
defined by $\tau$ with $\tau(u,a) = \Phi(u)(a)$ for all $u \in U,
a \in A.$ We write $ua$ instead of $u \otimes a \in H$ for all $u
\in U, a \in A.$

Let $\rho \in \widehat{\Lambda}$ and $\chi \in \widehat{\Gamma}$.
We recall that $L(\rho,\chi) = U_{\chi}/I(\rho,\chi)$ is a cyclic
left $H$-module and left $U$-module with generator
$$m=m_{\rho,\chi}=\text{ residue class of }1.$$
We denote the $H$-action on the quotient $L(\rho,\chi)$ by $\cdot$
and we write $um=u\cdot m$ for all $u \in U$. Recall that
$u\cdot_{\chi}u'=uu'$ for all $u,u' \in U$ by  \eqref{actionchi}.
Therefore for all $u \in U$, $um$ is the residue class of $u$ in
$L(\rho,\chi)$ and the $H$-action on $um$ is given for all $u' \in
U$ and $a \in A$ by
\begin{equation}\label{HactiononL}
u'\cdot um= u'um,\; a\cdot um = (a \cdot_{\chi}u)m.
\end{equation}
The latter holds since $a \cdot_{\chi}(u \cdot_{\chi}1) = a
\cdot_{\chi}u = (a \cdot_{\chi}u)\cdot_{\chi}1$.

We denote the braiding matrix of $W$ by $(q_{ij})$, where 
\begin{equation}\label{notation}
q_{ij} = \eta_j(z_i) \text{ for all } 1 \leq i,j \leq n.
\end{equation}

\begin{Prop}\label{actL}
Let $\rho \in \widehat{\Lambda}$ and $\chi \in \widehat{\Gamma}$.
Then $L(\rho,\chi)$ is the $k$-span of all
$$u_{i_1} \cdots u_{i_t} m,
\text{ where } 1 \leq i_1,\dots,i_t \leq n \text{ and } t \geq
0.$$ The $H$-action on these elements is given for all $z \in
\Lambda,g \in \Gamma,$ and  $1 \leq j \leq n$ by
\begin{enumerate}
\item [{\rm a)}] $z \cdot u_{i_1}\cdots u_{i_t} m\phantom{_j} =
(\eta_{i_1}\cdots \eta_{i_t} \rho)(z)u_{i_1} \cdots u_{i_t} m,$
\item [{\rm b)}] $g \cdot u_{i_1} \cdots u_{i_t} m \phantom{_j}=
(\chi_{i_1}^{-1}\cdots \chi_{i_t}^{-1} \chi)(g)u_{i_1} \cdots
u_{i_t} m,$ \item [{\rm c)}] $u_j \cdot u_{i_1} \cdots u_{i_t} m =
u_ju_{i_1}\cdots u_{i_t} m,$ \item [{\rm d)}] $a_j \cdot u_{i_1}
\cdots u_{i_t} m = \sum_{l=1}^t \alpha_l(j,i_1,\dots,i_t)u_{i_1}
\cdots u_{i_{l-1}} u_{i_{l+1}} \cdots u_{i_t}m,$
\end{enumerate}
with $\alpha_l(j,i_1,\dots,i_t)= \delta_{i_lj} l_j
\prod_{r=1}^{l-1}q_{i_{r}j}\left(1-\prod_{s=l+1}^{t}q_{i_{s}j}q_{ji_{s}}
\rho(z_j)\chi(g_j)\right).$
\end{Prop}
\pf We have seen in the proof of Theorem \ref{groupcase} that
$$z \cdot m = \rho(z)m.$$
Since $zu_i=\eta_i(z)u_iz$ for all $1 \leq i \leq n$, part a)
follows by the first equation of \eqref{HactiononL}. Note that
part a) implies that the elements $u_{i_1} \cdots u_{i_t} m$ span
$L(\rho,\chi)$.

Let $u = u_{i_1} \cdots u_{i_t}.$ Note that by \eqref{comultU}
\begin{align}
\Delta^2(u) &= u\sw1 \otimes u\sw2 \otimes u\sw3 \notag\\
&=(z_{i_1} \otimes z_{i_1} \otimes u_{i_1} + z_{i_1} \otimes u_{i_1} \otimes 1 + u_{i_1} \otimes 1 \otimes 1)\cdots\label{Delta2}\\
&\cdots(z_{i_t} \otimes z_{i_t} \otimes u_{i_t} + z_{i_t} \otimes
u_{i_t} \otimes 1 + u_{i_t} \otimes 1 \otimes 1).\notag
\end{align}

By \eqref{HactiononL} $g \cdot u m=  (g \cdot_{\chi} u) m,$  and
by \eqref{actionchi} and \eqref{Delta2}
\begin{align*}
g \cdot_{\chi} u &= \Phi(u\sw1)(g) u\sw2 \chi(g) \Phi(u\sw3)(g^{-1})\\
&= \varphi(z_{i_1} \cdots z_{i_t})(g) u_{i_1} \cdots u_{i_t}
\chi(g),
\end{align*}
since $\Phi(u\sw1)(g) = 0$ resp.  $\Phi(u\sw3)(g^{-1}) = 0$ if the
term $u\sw1$ resp. $u\sw3$ contains a factor $u_{i_l}.$ This
proves part b) since
$$\varphi(z_{i_1} \cdots z_{i_t})(g)= (\chi_{i_1}^{-1} \cdots \chi_{i_t}^{-1})(g)$$ by \eqref{conditionvarphi}.

Part c) is trivial, and to prove part d) we compute $a_j
\cdot_{\chi}u.$ Since $\Delta^2(a_j) = g_j \otimes g_j \otimes a_j
+ g_j \otimes a_j \otimes 1 + a_j \otimes 1 \otimes 1,$ and
$\chi(a_j)=0$ it follows from \eqref{actionchi} that
\begin{equation}\label{computation1}
a_j \cdot_{\chi}u = \Phi(u\sw1)(g_j) u\sw2 \chi(g_j)
\Phi(u\sw3)(-a_jg_j^{-1}) + \Phi(u\sw1)(a_j)u\sw2
\end{equation}
as $S^{-1}(a_j) = -a_jg_j^{-1}.$

To compute the first term in \eqref{computation1} we first note
that
\begin{equation}\label{computation2}
\Phi(u\sw1)(g_j)u\sw2 = q_{i_1j}\cdots q_{i_tj} u.
\end{equation}
For by part a) of Theorem \ref{groupcase} we have
$\Phi(u\sw1)(g_j) = 0$ if $u\sw1$ contains at least one factor
$u_i,$ hence $\Phi(u\sw1)(g_j)u\sw2 = \Phi(z_{i_1} \cdots
z_{i_t})(g_j)u.$ Next we see that
\begin{equation}\label{computation3}
u\sw1 \Phi(u\sw2)(a_jg_j^{-1}) = \sum_{l=1}^t u_{i_1} \cdots
u_{i_{l-1}}z_{i_l}u_{i_{l+1}}\cdots u_{i_t}
\Phi(u_{i_l})(a_jg_j^{-1})
\end{equation}
since $\Phi(u\sw2)(a_jg_j^{-1}) =0,$ if $u\sw2$ contains no or at
least two factors $u_{i_l}.$ Using \eqref{actionU},
\eqref{definitiondelta} and \eqref{notation} we obtain from
\eqref{computation3}
\begin{align}
&u\sw1 \Phi(u\sw2)(a_jg_j^{-1}) =\label{computation4}\\
&\sum_{l=1}^t \delta_{i_lj} l_jq_{jj}^{-1} q_{i_li_{l+1}} \cdots
q_{i_li_t} u_{i_1} \cdots u_{i_{l-1}} u_{i_{l+1}} \cdots u_{i_{t}}
z_{i_l}.\notag
\end{align}
We now use \eqref{computation4} in $\Delta$ applied to
\eqref{computation2} and obtain for the first term in
\eqref{computation1}
\begin{align}
&\Phi(u\sw1)(g_j) u\sw2 \chi(g_j) \Phi(u\sw3)(-a_jg_j^{-1})=\label{computation5}\\
- q_{i_1j}\cdots &q_{i_tj}\chi(g_j)\sum_{l=1}^t \delta_{i_lj} l_j
q_{jj}^{-1} q_{i_li_{l+1}} \cdots q_{i_li_t} u_{i_1} \cdots
u_{i_{l-1}} u_{i_{l+1}} \cdots u_{i_{t}} z_{i_l}.\notag
\end{align}
Similarly we compute the second term in \eqref{computation1}
\begin{equation}\label{computation6}
\Phi(u\sw1)(a_j)u\sw2 = \sum_{l=1}^t \delta_{i_lj} l_jq_{i_1i_l}
\dots q_{i_{l-1}i_l} u_{i_1} \cdots u_{i_{l-1}}u_{i_{l+1}} \cdots
u_{i_t}.
\end{equation}
Finally \eqref{computation1}, \eqref{computation5} and
\eqref{computation6} prove part d) in view of \eqref{HactiononL}
and part a).
 \epf  \medskip

The presentation of $L(\rho,\chi)$ in \eqref{presentationL} as the
subspace $U \cdot_{\rho}\chi$ of $A^*$ is helpful to actually
calculate the elements of $L(\rho,\chi)$ as linear functions on
$A.$ The next proposition in principle describes an inductive
procedure to calculate $L(\rho,\chi).$

\begin{Prop}\label{elementsofL}
Let $\rho \in \widehat{\Lambda}$ and $\chi \in \widehat{\Gamma}$,
and $1 \leq i,i_1,\dots,i_t \leq n,$ where $t\geq 1$. Define $u =
u_{i_1} \cdots u_{i_t}.$ Then
\begin{enumerate}
\item [{\rm a)}] $u_i \cdot_{\rho} \widetilde{\chi} =
(1-\rho(z_i)\chi(g_i)) \Phi(u_i)\widetilde{\chi}.$ \item [{\rm
b)}] $u_i \cdot_{\rho}(\Phi(u)\widetilde{\chi}) = \Phi(\left[u_iu
- \prod_{r=1}^{t}
q_{ii_{r}}\rho(z_i)\chi(g_i)uu_i\right])\widetilde{\chi}.$ \item
[{\rm c)}] $u_i^t \cdot_{\rho} \widetilde{\chi} =
\prod_{l=0}^{t-1}
\left(1-q_{ii}^l\rho(z_i)\chi(g_i)\right)\Phi(u_i^t)\widetilde{\chi}.$
\end{enumerate}
\end{Prop}
\pf Since $\Delta^2(u_i) = z_i \otimes z_i \otimes u_i + z_i
\otimes u_i \otimes 1 + u_i \otimes 1 \otimes 1,$ and $\rho(u_i) =
0,S(u_i)=-z_i^{-1}u_i$ it follows from \eqref{actionrho1} that for
any $p \in A^*$
\begin{equation}
u_i \cdot_{\rho}p = \Phi(z_i)\rho(z_i)p\Phi(-z_i^{-1}u_i) +
\Phi(u_i)p. 
\end{equation}
In particular
\begin{align}
u_i\cdot_{\rho}(\Phi(u)\widetilde{\chi}) &= \Phi(u_i)\Phi(u) \widetilde{\chi} - \rho(z_i)\Phi(z_i)\Phi(u)\widetilde{\chi}\Phi(z_i^{-1})\Phi(u_i)\label{el4}\\
&=\Phi(u_iu)\widetilde{\chi} - q_{ii_1} \cdots q_{ii_t}\rho(z_i)
\Phi(u) \widetilde{\chi} \Phi(u_i)\notag
\end{align}
since, as $\Alg(A,k)$ is abelian,
\begin{align*}
\Phi(z_i)\Phi(u) \widetilde{\chi} \Phi(z_i^{-1})  &= \Phi(z_i)\Phi(u)\Phi(z_i^{-1}) \widetilde{\chi}\\
&= \Phi(z_iuz_i^{-1}) \widetilde{\chi}\\
&=q_{ii_1} \cdots q_{ii_t} \Phi(u)\widetilde{\chi}.
\end{align*}

Since $\widetilde{\chi} \Phi(u_i)\widetilde{\chi}^{-1}$ and
$\chi(g_i) \Phi(u_i)$ are both
$(\varepsilon,\widetilde{\chi_i^{-1}})$-derivations taking the
same values on the generators $a_j, 1 \leq j \leq n,$ and $g \in
\Gamma$ of $A$ we see that
\begin{equation}\label{el5}
\widetilde{\chi} \Phi(u_i)\widetilde{\chi}^{-1} =\chi(g_i)
\Phi(u_i).
\end{equation}
Part b) now follows from \eqref{el4} and \eqref{el5}, and part a)
is the special case of part b) with $t=0.$

Part c) then follows by induction on $t.$ The case $t=1$ is part
a), and the induction step is
\begin{align*}
u_i^{t+1} \cdot_{\rho} \widetilde{\chi} &= u_i \cdot_{\rho}(u_i^t \cdot_{\rho}\widetilde{\chi})  \\
&=\prod_{l=0}^{t-1} (1-q_{ii}^l\rho(z_i)\chi(g_i))u_i\cdot_{\rho}(\Phi(u_i^t)\widetilde{\chi})\\
&= \prod_{l=0}^{t}
(1-q_{ii}^l\rho(z_i)\chi(g_i))\Phi(u_i^{t+1})\widetilde{\chi}
\end{align*}
since by part b) $u_i\cdot_{\rho}(\Phi(u_i^t)\widetilde{\chi})= (1
- q_{ii}^t\rho(z_i)\chi(g_i))\Phi(u_i^{t+1})\widetilde{\chi}.$
\epf

In the next Corollary we formulate a necessary condition for the
finiteness of the dimension of $L(\rho,\chi)$. We call a pair of
characters $(\rho,\chi) \in \widehat{\Lambda} \times \widehat{\G}$
{\em dominant} if there are natural numbers $m_i \geq 0,1 \leq i
\leq n,$ such that
\begin{equation}\label{dominant}
q_{ii}^{m_i} \rho(z_i) \chi(g_i) =1 \text{ for all } 1 \leq i \leq
n.
\end{equation}

\begin{Cor}\label{sufficient}
Let $\rho \in \widehat{\Lambda}$ and $\chi \in \widehat{\Gamma}$
and assume that the braiding of $W$ is generic. If $L(\rho,\chi)$
is finite-dimensional, then the pair $(\rho,\chi)$ is dominant.
\end{Cor}
\pf Let $1 \leq i \leq n.$ The elements $u_i^tm, t \geq 0,$ of
$L(\rho,\chi)$ are linearly dependent. By part a) of Proposition
\ref{actL}, the group $\Lambda$ acts on  $u_i^tm$ via the
character $\eta_i^t\rho.$ Since $q_{ii}$ is a root of unity,
$\eta_i^t\rho \neq \eta_i^{t'}\rho$ for all $t \neq t',$ there is
an integer $m_i \geq 0$ such that $u_i^tm =0$ for all $t >m_i,$
and $u_i^tm \neq 0$ for all $0 \leq t \leq m_i.$ It is well-known
that  $u_i^t \neq 0$ for all $t \geq 0$ since $q_{ii}$ is not a
root of unity (see for example \cite[Example 2.9]{ASSurvey}).
Since $\Phi$ is injective by Remark \ref{Phiinjective}, it follows
that $\Phi(u_i^t) \neq 0$ for all $t \geq 0$, and thus
$q_{ii}^{m_i} \eta(z_i)\chi(g_i) =1$ by part c) of Proposition
\ref{elementsofL} and \eqref{presentationL}.
 \epf  \medskip

As an example we consider the easiest case where $U$ and $A$ are
quantum linear spaces, that is \eqref{commutation1} holds.

\begin{Cor}\label{quantumlinear1}
Let $\rho \in \widehat{\Lambda}$ and $\chi \in \widehat{\Gamma}$
and assume that the braiding of $W$ is generic and satisfies
\begin{equation}\label{commutation1}
q_{ij}q_{ji}=1 \text{ for all }1 \leq i,j \leq n,i \neq j.
\end{equation}
Then for all $t_1, \dots,t_n \geq 0,$
$$(u_1^{t_1}\cdots u_n^{t_n}) \cdot_{\rho} \widetilde{\chi} = \prod_{i=1}^n \prod_{l_i = 0}^{t_i -1} (1 - q_{ii}^{l_i} \rho(z_i)\chi(g_i)) \Phi(u_1^{t_1}\cdots u_n^{t_n})\widetilde{\chi}.$$
\end{Cor}
\pf It follows from \eqref{commutation1} that
\begin{equation}\label{commutation2}
u_iu_j = q_{ij}u_ju_i\text{ for all }1 \leq i,j \leq n,i \neq j,
\end{equation}
since the skew-commutator $u_iu_j - q_{ij}u_ju_i$ is primitive in
$\mathfrak{B}(W).$

We prove by induction on $1\leq j \leq n$ that
$$(u_j^{t_j}\cdots u_n^{t_n}) \cdot_{\rho} \widetilde{\chi} = \prod_{i=j}^n \prod_{l_i = 0}^{t_i -1} (1 - q_{ii}^{l_i} \rho(z_i)\chi(g_i)) \Phi(u_j^{t_j}\cdots u_n^{t_n})\widetilde{\chi}$$
for all $t_j,\dots,t_n \geq 0.$

The case $j = n$ is part c) of Proposition \ref{elementsofL}.
Assume the formula for $j+1 \leq n.$ Then we prove the claim for
$j$ by induction on $t_j$. If $t_j=0,$ then the formula is true by
assumption, and for all $t_j \geq 0$ we see by induction on $t_j$
that
\begin{align}
(u_j^{t_j+1} \cdots u_n^{t_n}) \cdot_{\rho}\widetilde{\chi}
&=u_j \cdot_{\rho}\left((u_j^{t_j}\cdots u_n^{t_n}) \cdot_{\rho} \widetilde{\chi}\right)\notag\\
&=\prod_{i=j}^n \prod_{l_i = 0}^{t_i -1} (1 - q_{ii}^{l_i}
\rho(z_i)\chi(g_i)) u_j \cdot_{\rho}(\Phi(u_j^{t_j} \cdots
u_n^{t_n})\widetilde{\chi}).\label{induction}
\end{align}
By part b) of Proposition \ref{elementsofL},
\begin{eqnarray*}
\lefteqn {u_j\cdot_{\rho}\left(\Phi(u_j^{t_j} \cdots u_n^{t_n})\widetilde{\chi}\right)} \\
&=&\left(\Phi(u_j^{t_j+1} \cdots u_n^{t_n}) - q_{jj}^{t_j} q_{j,j+1}^{t_{j+1}} \cdots q_{jn}^{t_n} \rho(z_j)\chi(g_j) \Phi(u_j^{t_j} \cdots u_n^{t_n}u_j)\right)\widetilde{\chi}\\
&=& \left(1 - q_{jj}^{t_j}
\rho(z_j)\chi(g_j)\right)\Phi(u_j^{t_j+1} \cdots
u_n^{t_n})\widetilde{\chi},
\end{eqnarray*}
since $u_j^{t_j} \cdots u_n^{t_n} u_j = q_{j+1,j}^{t_{j+1}} \cdots
q_{nj}^{t_n} u_j^{t_j+1} \cdots u_n^{t_n}$ by
\eqref{commutation2}, and all the $q_{jl}$ and $q_{lj}$ except for
$l=j$ cancel by \eqref{commutation1}. Hence the claim follows from
\eqref{induction}. \epf

\begin{Cor}\label{quantumlinear2}
Let $\rho \in \widehat{\Lambda}$ and $\chi \in \widehat{\Gamma}$
and assume that the braiding of $W$ is generic and satisfies
\eqref{commutation1}. Then $L(\rho,\chi)$ is finite-dimensional if
and only if the pair $(\rho,\chi)$ is dominant. In this case, the
elements $u_1^{t_1} \cdots u_n^{t_n}m, 0 \leq t_i \leq m_i$ for
all $1 \leq i \leq n,$ form a basis of $L(\rho,\chi).$
\end{Cor}
\pf It is well-known that by \eqref{generic} and
\eqref{commutation1} $\mathfrak{B}(W)$ is generated by
$u_1,\dots,u_n$ with defining relations 
$$u_iu_j - q_{ij}u_ju_i =0
\text{ for all } 1 \leq i,j \leq n, i \neq j,$$ 
and has PBW-basis
$u_1^{t_1} \cdots u_n^{t_n}, t_1,\dots,t_n \geq 0$ (see for
example \cite[Theorem 4.2]{AS3}).

If $L(\rho,\chi)$ is finite-dimensional, then $(\rho,\chi)$ is
dominant by Corollary \ref{sufficient}. Conversely, assume that
$(\rho,\chi)$ is dominant. Since $\Phi$ is injective by Remark
\ref{Phiinjective}, and the elements $u_1^{t_1} \cdots u_n^{t_n},
0 \leq t_i \leq m_i$ for all $1 \leq i \leq n,$ of the PBW-basis
are linearly independent it follows from Corollary
\ref{quantumlinear1} that the elements $u_1^{t_1} \cdots
u_n^{t_n}m, 0 \leq t_i \leq m_i$ for all $1 \leq i \leq n,$ are
linearly independent, hence a basis of $L(\rho,\chi)$ since
$u_1^{t_1} \cdots u_n^{t_n} = 0$  if $t_i > m_i$ for one $i$. \epf



We now extend the finiteness criterion for quantum linear spaces
in Corollary \ref{quantumlinear2} to braidings  of finite Cartan
type.

\begin{Theorem}\label{dimfinite}
Assume that the braiding of $W$ is generic and of finite Cartan
type. let $e_1,\dots,e_n >1$ be natural numbers. Then
$$\mathfrak{B}(W)/(\sum_{1 \leq i \leq n}\mathfrak{B}(W)u_i^{e_i})$$
is finite-dimensional.
\end{Theorem}
\pf As explained in Section \ref{Sectiongeneral} there are $d_i
\in \{1,2,3\},1 \leq i \leq n$, and $q_{J} \in k, J \in
\mathcal{X}$, such that 
$$q_{ii}=q_J^{2d_i}, d_ia_{ij} =d_ja_{ji}$$
for all $J \in \mathcal{X},i,j \in J$. We define $(q'_{ij})$ by \eqref{q'}.

Let $\mathbb{Z}[I]$ be the free abelian group of rank $n$ with
basis $\alpha_1,\dots,\alpha_n$.  Define characters $\psi,\psi'
\in \widehat{\mathbb{Z}[I]}$ by
\begin{equation}
\psi_j(\alpha_i)=q_{ij},\psi'_j(\alpha_i) = q'_{ij} \text{ for all
} 1 \leq i,j \leq n.
\end{equation}
Let $X, X' \in {\YDI}$ be Yetter-Drinfeld modules with bases $x_i \in X^{\psi}_{\alpha_i}$ and $x'_i
\in {X'}^{\psi'_i}_{\alpha_i},1 \leq i \leq n$.

Then
\begin{equation}
\mathfrak{B}(W) \to \mathfrak{B}(X), u_i \mapsto x_i,1 \leq i \leq
n,
\end{equation}
defines an isomorphism of algebras and coalgebras since the
braiding matrices of $W$ and $X$ coincide. Hence it suffices to
prove the theorem for $X$ instead of $W$.

By \eqref{twist} and \cite[Proposition 2]{AS3} there is a
2-cocycle 
$$\sigma : \mathbb{Z}[I] \times \mathbb{Z}[I] \to
k^{\times}$$ of the group $\mathbb{Z}[I]$ and a $k$-linear
isomorphism
\begin{align}
&\Psi :  \mathfrak{B}(X) \to \mathfrak{B}(X') \text{ with }\Psi(x_i)=x'_i,1 \leq i \leq n\label{Psi}\\
\intertext{such that for all $\alpha, \beta \in \mathbb{Z}[I]$ and
$x \in \mathfrak{B}(X)_{\alpha},y \in \mathfrak{B}(X)_{\beta}$,}
&\Psi(xy)=\sigma(\alpha,\beta) \Psi(x) \Psi(y).\label{Psisigma}
\end{align}
Let $\mathcal{X} = \{I_1,\dots,I_t\}$ be the set of connected
components of $\{1,2, \dots, n\}$. For all $1 \leq s \leq t$
we let $X_s= \oplus_{i \in I_s} k x_i$. Then the natural map
\begin{align*}
\mathfrak{B}(X_1)\otimes \cdots \otimes \mathfrak{B}(X_t)
 \to \mathfrak{B}(X)
\end{align*}
given by multiplication is an isomorphism by \cite[Lemma
1.4]{AS3}. Since $x_ix_j = q_{ij}x_jx_i$ for all $1 \leq i,j \leq n, i \nsim j$, it induces a surjective linear map from
$$\mathfrak{B}(X_1)/(\sum_{i \in I_1} \mathfrak{B}(X_1)x_i^{e_i}) \otimes \cdots \otimes
\mathfrak{B}(X_t)/(\sum_{i \in I_t} \mathfrak{B}(X_t)x_i^{e_i})$$
to
 $$\mathfrak{B}(X)/(\sum_{1 \leq i \leq n}\mathfrak{B}(X)x_i^{e_i}).$$
Thus we have reduced the claim to the connected case.

If the Dynkin diagram of $(a_{ij})$ is connected, let $q =q_{I}$.
Then by results of Lusztig \cite[Section 37]{L} and Rosso
\cite[Theorem 15]{R} $$\mathfrak{B}(X')=k \langle x'_1,\dots,x'_n
\mid (\ad_{c'}(x'_i))^{1-a_{ij}}(x'_j)=0, i \neq j\rangle$$ is
isomorphic to $U_q^-(\mathfrak{g})$, where $\mathfrak{g}$ is the
semisimple Lie algebra with Cartan matrix $(a_{ij})$. The elements
$x'_i$ correspond to the $F_i$. Hence by \cite[Proposition 6.3.4]{L} or \cite[Proposition 5.9]{J}
$$\mathfrak{B}(X')/(\sum_{1 \leq i \leq n}\mathfrak{B}(X') {x'}_i^{e_i})$$
is finite-dimensional. By \eqref{Psi} and \eqref{Psisigma}, $\Psi$
induces an isomorphism of vector spaces
$$\mathfrak{B}(X)/(\sum_{1 \leq i \leq n}\mathfrak{B}(X) {x}_i^{e_i})\cong \mathfrak{B}(X')/(\sum_{1 \leq i \leq n}\mathfrak{B}(X') {x'}_i^{e_i}),$$
and the theorem is proved. \epf

\medskip

We can now prove the main result of this section.

\begin{Theorem}\label{maingeneric}
Assume that the braiding of $W$ is generic and of finite Cartan
type. Then the map
$$\{(\rho,\chi) \in \widehat{\Lambda} \times \widehat{\G} \mid (\rho,\chi) \text{ dominant}\} \to \Irr((U \otimes A)^{\sigma}),$$
given by $(\rho,\chi) \mapsto [L(\rho,\chi)]$, is bijective.
\end{Theorem}
\pf Let $(\rho,\chi) \in \widehat{\Lambda} \times \widehat{\G}$ be
a dominant pair. By Theorem \ref{dim1}, Theorem \ref{modules} and
Proposition \ref{sufficient} it suffices to show that
$L(\rho,\chi)$ is finite-dimensional.

By definition $q_{ii}^{m_i} \rho(z_i) \chi(g_i) =1$ for all $1
\leq i \leq n$, where $m_i \geq 0$ are natural numbers. Then it
follows from Proposition \ref{elementsofL} c) that $u_i^{m_i+1}
\in I(\rho,\chi)$ for all $1 \leq i \leq n$. The natural map
$$\mathfrak{B}(W) \subset U \to L(\rho,\chi) = U/I(\rho,\chi)$$
is surjective, since $z \cdot m = \rho(z)m$ for all $z \in
\Lambda$ by Proposition \ref{actL} a). Thus we have a surjective
linear map
$$\mathfrak{B}(W)/(\sum_{1 \leq i \leq n} \mathfrak{B}(W)u_i^{m_i+1}) \to L(\rho,\chi),$$
 and $L(\rho,\chi)$ is finite-dimensional by Theorem \ref{dimfinite}.
\epf

\section{Application to pointed Hopf algebras given by a datum of finite Cartan type and linking parameters}\label{SectionCartan}

Let $\mathcal{D} = \mathcal{D}(\Gamma, (g_i)_{1 \leq i \leq
\theta}, (\chi_i)_{1 \leq i \leq \theta}, (a_{ij})_{1 \leq i,j
\leq \theta})$ be a datum of finite Cartan type and suppose that
$\lambda$ is a family of linking parameters for $\mathcal{D}$ in
the sense of \cite{AS2}. We apply the theory developed in Sections \ref{Sectiongeneral} and \ref{computations} to the infinite-dimensional Hopf algebra
$U(\mathcal{D},\lambda)$ in general and to the finite-dimensional
version $u(\mathcal{D},\lambda)$ under some restrictions. See definitions \ref{DefU} and \ref{Defu}.

Each is a quotient of an $H=(U \otimes A)^{\sigma}$ where the
2-cocycle $\sigma$ is determined by the linking parameters. The
finite-dimensional irreducible left modules for
$U(\mathcal{D},\lambda)$ are parameterized by a subset of
$\widehat{\G}$ and the irreducible modules for
$u(\mathcal{D},\lambda)$ by $\widehat{\G}$; these irreducible
modules arise from $L(\rho,\chi)$'s defined for $H$ where $\chi$
determines $\rho$.

\subsection{The Linking Graph}\label{Linking}

Let
$$\mathcal{D} = \mathcal{D}(\Gamma, (g_i)_{1 \leq i \leq \theta}, (\chi_i)_{1 \leq i \leq \theta}, (a_{ij})_{1 \leq i,j \leq \theta})$$
be a  of finite Cartan type. Note that by \eqref{Cartan}
\begin{equation}\label{Cartan1}
q_{ii}^{a_{ij}}=q_{jj}^{a_{ji}} \text{ for all } 1\leq i,j \leq
\theta.
\end{equation}

Let $1 \leq i,j \leq \theta$. Then $i\sim j$, that is, $i$ and $j$
lie in the same connected component of the Dynkin diagram of
$(a_{ij})$ if and only if $i=j$ or there are indices $1 \leq i_1,
\dots,i_k \leq \theta, i_l \neq i_{l+1} \text{ for all } 1 \leq l
<k$, such that $a_{i_li_{l+1}} \neq 0$ for all $1 \leq l <k$. In
this case we write
\begin{align}
a(i,j) &= a_{i_1i_2} a_{i_2i_3} \cdots a_{i_{k-1}i_k}, b(i,j) =a_{i_2i_1} a_{i_3i_2} \cdots a_{i_{k}i_{k-1}},\label{ab}\\
\intertext{and by \eqref{Cartan1} we have}
q_{ii}^{a(i,j)}&=q_{jj}^{b(i,j)}.\label{ab1}
\end{align}
Note that $a(i,j)$ and $b(i,j)$ depend on the choice of the
sequence $i_1,\dots,i_k$.

A family $ \lambda = (\lambda_{ij})_{1 \le i, j \le \theta, \,
i\not\sim j}$ of elements in $k$ is called a {\em family of
linking parameters for $\D$} if the following conditions are
satisfied for all $1 \le i, j \le \theta, \, i\not\sim j$:
\begin{align}
\text{ If }g_ig_j &=1  \text{ or }\chi_i \chi_j \neq \varepsilon  \text{ then } \lambda_{ij} =0;\label{Deflinkable}\\
\lambda_{ji} &= - q_{ji} \lambda_{ij}.\label{linkingnumbering}
\end{align}
This definition is a formal extension of \cite[Section 5.1]{AS2}
where $\lambda_{ij}$ was only defined for $i<j$ and $i \not\sim
j.$ Note that by \eqref{linkingnumbering} 
$$\lambda_{ij} = -
q_{ji}^{-1}\lambda_{ji} = -q_{ij} \lambda_{ji}$$ since
$q_{ij}q_{ji} = 1$ for all $i \not\sim j.$ Note that
\eqref{Deflinkable} and \eqref{linkingnumbering} are met when
$\lambda_{ij}=0$ for all $i,j$. We let $0$ denote this family of
linking parameters.

Vertices $1\leq i,j \leq \theta$ are called {\em linkable} if $i
\not\sim j,$ $g_ig_j \neq 1$ and $\chi_i \chi_j = \varepsilon$.
Then (see \cite[Section 5.1]{AS2})
\begin{equation}\label{linkable}
q_{ii} = q_{jj}^{-1} \text{ if }i,j \text{ are linkable}.
\end{equation}

The next lemma is \cite[Lemma 5.6]{AS2}.
\begin{Lem}\label{linkingLemma}
Let $\mathcal{D} = \mathcal{D}(\Gamma, (g_i)_{1 \leq i \leq
\theta}, (\chi_i)_{1 \leq i \leq \theta}, (a_{ij})_{1 \leq i,j
\leq \theta})$ be a datum of finite Cartan type, and assume that
$\ord(q_t) >3$ for all $1 \leq t \leq \theta.$
\begin{enumerate}
\item \phantom{a} If vertices $i,k$ and $j,l$ are linkable, then
$a_{ij}=a_{kl}.$\label{Cartanlinkable} 
\item \phantom{a} A vertex $i$ cannot be linkable to two different vertices
$j,k.$\label{unique}
\end{enumerate}
\end{Lem}

\medskip

\noindent We say that linkable vertices $1\leq i,j \leq \theta$
are {\em linked} if $\lambda_{ij} \neq 0.$

Let $\lambda$ be a family of linking parameters for a datum of
Cartan type $\mathcal{D} = \mathcal{D}(\Gamma, (g_i)_{1 \leq i
\leq \theta}, (\chi_i)_{1 \leq i \leq \theta}, (a_{ij})_{1 \leq
i,j \leq \theta}).$ We define the {\em linking graph} of $(\D,
\lambda)$ as follows: The set of its vertices is $\mathcal{X}$,
the set of all connected components of $I=\{1,\dots,n\}$. There is
an edge between $J_1,J_2 \in \mathcal{X}$ if and only there are
elements $i \in J_1,j\in J_2$ such that $i,j$ are linked. Recall
that a graph is called {\em bipartite} if the set of its vertices
is the disjoint union of subsets $X^+,X^-$ such that there is no
edge between vertices in $X^+$ or in $X^-$.

\begin{Lem}\label{linkinggraph}
Let $\mathcal{D} = \mathcal{D}(\Gamma, (g_i)_{1 \leq i \leq
\theta}, (\chi_i)_{1 \leq i \leq \theta}, (a_{ij})_{1 \leq i,j
\leq \theta})$ be a datum of Cartan type, and $\lambda$ a family
of linking parameters for $\D.$ Assume one of the following
conditions:
\begin{enumerate}
\item \phantom{a} The Cartan matrix $(a_{ij})$ is simply laced,
that is $a_{ij} \in \{0,-1\}$ for all $1 \leq i,j \leq \theta,i
\neq j,$ and $\ord(q_{ii})$ is odd for all $1 \leq i \leq \theta.$
\item \phantom{a} For all $1 \leq i \leq \theta,$ $q_{ii}$ is not a
root of unity.
\end{enumerate}
Then the linking graph of $(\D, \lambda)$ is bipartite.
\end{Lem}
\pf A graph is bipartite if and only if it contains no cycle of
odd length \cite[Proposition 0.6.1]{D}. Assume the linking graph
has a cycle of odd length n. Then there are connected components
$J_1,\dots,J_n,$ and $i_l \in J_l, 1 \leq l \leq n,$ and $j_l \in
J_{l}, 2 \leq l \leq n+1,J_{n+1}=J_1$ with $\lambda_{i_lj_{l+1}}
\neq 0$  for all $1 \leq l \leq n$.  Let $i_{n+1} =i_1$. For
simplicity we write $q_i=q_{ii}$ for all $i$. Since $i_l$ and
$j_l$ are in the same connected component for all $2 \leq l \leq
n+1$, by \eqref{ab1} $q_{i_l}^{a_l} = q_{j_l}^{b_l}$ for all $2
\leq l \leq n+1,$ where $a_{l} = a(i_l,j_l)$ and $b_l =
b(i_l,j_l)$ are chosen as in \eqref{ab}. By \eqref{linkable},
$q_{i_l} = q_{j_{l+1}}^{-1}$ for all $1 \leq l \leq n.$ Hence
$q_{i_l}^{a_l} = q_{i_{l-1}}^{-b_l}$ for all $2 \leq l \leq n+1$,
and
$$q_{i_{n+1}}^{a_{n+1}  \cdots a_2} = q_{i_1}^{(-1)^n b_{n+1} \cdots b_2}.$$
Since $i_{n+1} =i_1,$ and $n$ is odd, the order of $q_{i_1}$ must
divide 
$$a_{n+1} \cdots a_2 + b_{n+1} \cdots b_2.$$ 
This is
impossible in both cases (1) and (2).
 \epf

 \medskip

However, in general the linking graph of $(\D,\lambda)$ is not
necessarily bipartite.
\begin{Ex}\label{B2}
One can link an odd number of copies of $B_2$ in a circle, where
the group $\Gamma$ is $\mathbb{Z}^{2n}$ or $(\mathbb{Z}/(N))^{2n}$
such that $N$ is an integer dividing $1+2^n,$ and where
$g_1,\dots,g_{2n}$ are the canonical basis elements of $\Gamma.$
\epf\end{Ex}

\subsection{The Infinite-dimensional Case}\label{infinite}

In this section we fix a generic Cartan datum of finite type
$\mathcal{D} = \mathcal{D}(\Gamma, (g_i)_{1 \leq i \leq \theta},
(\chi_i)_{1 \leq i \leq \theta}, (a_{ij})_{1 \leq i,j \leq
\theta}),$ and a family $ \lambda = (\lambda_{ij})_{1 \le i, j \le
\theta, \, i\not\sim j}$ of linking parameters for $\D.$ We define
$q_{ij} = \chi_j(g_i)$for all $1 \leq i,j \leq \theta$.

Let $X\in {\YDG}$ be the vector space with basis $x_i \in
X_{g_i}^{\chi_i}, 1 \leq i \leq \theta.$ The tensor algebra $T(X)$
is an algebra in $\YDG$. Thus the smash product $T(X) \# k[\G]$ is
a biproduct and thus has a Hopf algebra structure. We identify
$T(X)$ with the free algebra $k\langle x_1, \dots,x_{\theta}
\rangle.$

\begin{Def}\label{DefU}
Let $U(\D,\lambda)$ be the quotient Hopf algebra of the smash
product $k\langle x_1,\dots,x_{\theta}\rangle \# k[\Gamma]$ modulo
the ideal generated by
\begin{align}
&\ad_c(x_i)^{1-a_{ij}}(x_j), \text{ for all } 1\leq i,j \leq \theta, i \sim j, i \neq j,&\label{Serrerelations}\\
&x_i x_j - q_{ij} x_j x_i - \lambda_{ij}(1 - g_i g_j), \text{ for
all } 1 \le i < j \le \theta, \, i\not\sim
j.&\label{linkingrelations}
\end{align}
\end{Def}
\noindent Here, the $k$-linear endomorphism $\ad_c(x_i)$ of the
free algebra is given for all $y$ by the braided commutator
$\ad_c(x_i)(y)=x_iy - (g_i \cdot y)x_i$.

We denote the images of $x_i$ and $g\in \Gamma$ in $U(\D,\lambda)$
again by $x_i$ and $g$. Note that in $U(\D,\lambda)$
$$x_i x_j - q_{ij} x_j x_i = \lambda_{ij}(1 - g_i g_j), \text{ for all } 1 \le i,j \le \theta, \, i\not\sim j$$
by \eqref{linkingnumbering}.

The elements in \eqref{Serrerelations} and
\eqref{linkingrelations} are skew-primitive. Hence $U(\D,\lambda)$
is a Hopf algebra with
$$\Delta(x_i) = g_i \o x_i + x_i \o 1,\; 1 \leq i \leq \theta.$$
Then $U(\D,0)$ is the biproduct $\mathfrak{B}(X) \# k[\G]$.

In the next Lemma we write down the relations
\eqref{Serrerelations} explicitly, and we give an equivalent
formulation replacing the elements $x_i$ by $x_ig_i^{-1}$. The
elements $x_i$ correspond to the usual $E_i$, and the elements
$x_ig_i^{-1}$ to the $F_i$ in $U_q(\mathfrak{g})$, $\mathfrak{g}$
a semisimple Lie algebra. We recall the box notation for
$q$-integers (see \cite[1.3]{L}, \cite[Chapter 0]{J}). Let $v$ be
an indeterminate. For all integers $n \geq 0$ let $$[n]= \frac{v^n
- v^{-n}}{v-v^{-1}},$$ $[0]^! =1,$ and $[n]^!=[1][2]\cdots[n]$ if
$n\geq 1$. For all $0 \leq i \leq n$ define
$$\begin{bmatrix}n\\i\end{bmatrix}=\frac{[n]^!}{[i]^![n-i]^!}.$$
If $ 0\neq q \in k$, then $[n]_q$ and
$\begin{bmatrix}n\\i\end{bmatrix}_q$ are the specializations at $v
\mapsto q$.
\begin{Lem}\label{equivalentSerre}
Let $J \in \mathcal{X}$ be a connected component of
$\{1,2,\dots,\theta\}$, and define $F_i=x_ig_i^{-1}$ for all $i
\in J$. Choose $q_J \in k$ such that \eqref{DJ} is satisfied. Then
the following equalities hold in the smash product $k\langle
x_1,\dots,x_{\theta}\rangle \# k[\Gamma]$ for all $i,j \in J$ and
any natural number $a$:
\begin{align}
(\ad_cx_i)^a(x_j) =& \sum_{s=0}^a (-q_J^{d_i(a-1)}q_{ij})^s
\begin{bmatrix}
a\label{SerreEversion}\\
s
\end{bmatrix}_{q_J^{d_i}}
x_i^{a-s}x_jx_i^s,\\
(\ad_cx_i)^a(x_j) =& \sum_{s=0}^a (-q_J^{d_i(a-1)}q_{ji})^s
\begin{bmatrix}
a\\
s
\end{bmatrix}_{q_J^{d_i}}
F_i^{a-s}F_jF_i^s\label{SerreFversion}\\
&\cdot g_i^ag_j(q_J^{d_i(a-1)}q_{ij})^a.\notag
\end{align}
\end{Lem}
\pf We extend $\ad_cx_i$ to a map on the smash product by
$$\ad_cx_i=L_{x_i}-R_{x_i}\sigma_i,$$ where $L_{x_i}$ and
$R_{x_i}$ denote left and right multiplication with $x_i$, and
$\sigma_i$ is the inner automorphism defined for all $ y \in k\langle
x_1,\dots,x_{\theta}\rangle \# k[\Gamma]$ by $\sigma_i(y) =
g_iy g_i^{-1}$. Since
$(R_{x_i}\sigma_i) L_{x_i}= q_{J}^{2d_i}
L_{x_i}(R_{x_i}\sigma_i)$, we can compute
$(L_{x_i}-R_{x_i}\sigma_i)^a$ by the $q$-binomial formula
\cite[1.3.5]{L}, and \eqref{SerreEversion} follows from the
equality $(R_{x_i}\sigma_i)^n=R_{x_i^n}\sigma_i^n q_J^{d_in(n-1)}$
for all natural numbers $n$. Then \eqref{SerreFversion} follows
from \eqref{SerreEversion} by a somewhat tedious computation where
we write $x_i^{a-s}x_jx_i^s$ as a multiple of $F_i^{a-s}F_jF_i^s$.
\epf

Suppose that $\lambda \neq 0.$ In \cite{AS3} a glueing process was
used to build  $U(\D,\lambda)$ inductively by adding one connected
component at a time and thus to obtain a basis of the algebra. The
methods of the present paper to parameterize the
finite-dimensional irreducible modules do not apply to this
description of $U(\D,\lambda)$.

Since by Lemma \ref{linkinggraph} the linking graph of
$(\D,\lambda)$ is bipartite, we can give another description of
$U(\D,\lambda)$ by one glueing only and show that the algebra is a
quotient of $(U \otimes A)^{\sigma}$ where some central
group-likes are identified with 1 and $U$ and $A$ are biproducts.
As a consequence Theorem \ref{maingeneric} applies.

In the inductive construction of \cite{AS3} each stage yields a
quotient of the form just described; however, $U$ is not a
biproduct and the finite-dimensional irreducible $U$-modules are
not necessarily one-dimensional.

We now give our construction of $U(\D,\lambda)$. By Lemma
\ref{linkingLemma} there are non-empty disjoint subsets
$\mathcal{X}^-, \mathcal{X}^+\subseteq \mathcal{X}$ with
$\mathcal{X} = \mathcal{X}^- \cup \mathcal{X}^+,$ an integer $n
\geq 1$ and an injective map $t : \{1,\dots,2n\} \to
\{1,\dots,\theta\},$ such that
$$t(i) \in I^-:=\bigcup_{J \in \mathcal{X}^-} J, \;t(n+i)\in I^+:=\bigcup_{J \in \mathcal{X}^+} J,$$ for all $1 \leq i \leq n,$ and such that $(t(i),t(n+i))$ and $(t(n+i),t(i))$, where $1 \leq i \leq n$ are all the linked pairs of elements in $\{1,\dots,\theta\}.$

Let $\Lambda$ be the free abelian group with basis $z_i,i \in
I^-.$ For all $j \in I^-$ let $\eta_j \in \widehat{\Lambda}$ be
defined by $\eta_j(z_i) = \chi_j(g_i)$ for all $i \in I^-.$

Let $W \in {\YDL}$ with basis $u_i \in W_{z_i}^{\eta_i}, i \in
I^-,$ and $V \in {\YDG}$ with basis $a_j \in V_{g_j}^{\chi_j}, j
\in I^+.$ We define
$$U=\mathfrak{B}(W) \# k[\Lambda] \text{ and }A=\mathfrak{B}(V) \# k[\G].$$

\begin{Theorem}\label{twistingbiproduct}
Assume the situation above. Then for all $i \in I^-$ there is a unique algebra map
\begin{equation}
\gamma_i : A \to k \text{ with }\gamma_i(a_j) = 0, \gamma_i(g) =
\chi_i(g)
\end{equation}
for all $j \in I^+$ and $g \in \Gamma,$ and a unique
$(\varepsilon,\gamma_i)$-derivation
\begin{equation}
\delta_i : A \to k \text{ with } \delta_i(a_j) = \lambda_{ji} ,
\delta_i(g)  = 0
\end{equation}
for all $j \in I^+$ and $g \in \Gamma.$ Moreover there is a Hopf
algebra map $\Phi : U \to  A^{o\,\cop}$ determined by
\begin{equation}
\Phi(z_i) = \gamma_i, \;\Phi(u_i) = \delta_i
\end{equation}
for all $i \in I^-.$

Let $\sigma$ be the 2-cocycle corresponding to $\Phi$ (see Section
1). The group-like elements $z_i \otimes g_i^{-1}, i \in I^-,$ are
central in $(U \otimes A)^{\sigma},$ and there is an isomorphism
of Hopf algebras
\begin{equation}\label{IsoD}
\Psi : U(\D,\lambda) \xrightarrow{\cong} (U \otimes
A)^{\sigma}/(z_i \otimes g_i^{-1} - 1 \otimes 1 \mid i \in I^-)
\end{equation}
mapping $x_i$ with $i \in I^-, x_j$ with $j\in I^+,$ respectively
$g \in \Gamma$ onto the residue classes of $u_i \otimes 1,1
\otimes a_j$, respectively $1 \otimes g.$
\end{Theorem}
\pf It can be checked directly as in \cite[Lemma 5.19]{AS2} that
the maps $\gamma_i, \delta_i$ and $\Phi$ are well-defined by
working with the defining relations.

To see that $\Phi$ is well-defined without checking the relations
we alternatively can apply Theorem \ref{groupcase}. We define
$\varphi : \Lambda \to \widehat{\Gamma}, (l_i)_{i\in I^-},$ and $s
: I^- \to I^+$ for all $i \in I^-$ by
\begin{align*}
\varphi(z_i) &= \chi_i,\\
l_i&=\begin{cases}
\lambda_{t(n+k),t(k)}& \text{if $i=t(k)$ for some $1 \leq k \leq n$},\\
0& \text{otherwise},
\end{cases}\\
s(i)&=\begin{cases}
t(n+k)& \text{if $i=t(k)$ for some $1 \leq k \leq n$},\\
j_0& \text{otherwise},
\end{cases}
\end{align*}

\noindent where $j_0$ is any element in $I^+.$

Hence $l_i \delta_{s(i),j} = \lambda_{ji}$ for all $i \in I^-,j
\in I^+.$ We have to check the conditions in \eqref{EqVarPhI2}.
Let $i \in I^-$ with $l_i \neq 0,$ that is $i=t(k)$ for some $1
\leq k \leq n.$ Then $\varphi(z_{t(k)}) = \chi_{t(k)} =
\chi_{t(n+k)}^{-1}$ since $t(k)$ and $t(n+k)$ are linked. This
proves the first part of \eqref{EqVarPhI2} since $i=t(k)$ and
$s(i) = t(n+k).$ The second part of \eqref{EqVarPhI2} says that
$\eta_i(z_j) =\varphi(z_j)(g_{s(i)}) = \chi_j(g_{s(i)})$ for all
$j \in I^-$. Since $j$ and $s(i)$ are in different connected
components, it follows from the Cartan condition that
$\chi_j(g_{s(i)}) = \chi_{s(i)}(g_j)^{-1}.$ Since $i$ and $s(i)$
are linked, $\chi_i = \chi_{s(i)}^{-1}.$ Hence $\chi_j(g_{s(i)}) =
\chi_i(g_j) = \eta_i(z_j)$ by definition of $\eta_i.$

The remaining claims of the theorem follow by direct calculations
as in \cite[Theorem 5.17, end of the proof]{AS2}.
 \epf  \medskip

We continue with the situation of Theorem \ref{twistingbiproduct}.
Let $\rho \in \widehat{\Lambda}, \chi \in \widehat{\G}.$ By
\cite[Lemma 4.3]{Pairs} the left $(U \otimes A)^{\sigma}$-module
$L(\rho,\chi)$ is annihilated by all $z_i \otimes g_i^{-1} - 1
\otimes 1, i \in I^-$, if and only if
\begin{equation}\label{defrho}
\rho(z_i)=\chi(g_i) \text{ for all }i \in I^-.
\end{equation}
Thus for any $\chi \in \widehat{\G}$ define  $\rho \in
\widehat{\Lambda}$ by \eqref{defrho}, and define
$$L(\chi)=L(\rho,\chi)$$
as a left module over $U(\D,\lambda)$ using the isomorphism
$$\Psi : U(\D,\lambda) \xrightarrow{\cong} (U \otimes A)^{\sigma}/(z_i \otimes g_i^{-1} - 1 \otimes 1 \mid  i \in I^-).$$

We call a character $\chi \in\widehat{\G}$ {\em dominant} (for
$\D$, $\lambda$ and $I^+$) if for all pairs $(i,j)$ of linked
elements $1 \leq i,j \leq \theta,$ where $i \in I^-,j \in I^+,$
there are integers $m_i \geq 0$ such that
\begin{equation}\label{dominantchi}
q_{ii}^{m_i} \chi(g_ig_j) =1.
\end{equation}
Note that this definition depends on the choice of the sets $I^-$
and $I^+$, since $q_{ii} = q_{jj}^{-1}$ for linked pairs $(i,j)$.

\begin{Theorem}\label{mainD}
Let $\D$ be a generic datum of finite Cartan type, and $\lambda$ a
family of linking parameters for $\D$.
\begin{enumerate}
\item The map
$$\{\chi \in \widehat{\G} \mid \chi \text{ dominant }\} \to \Irr(U(D,\lambda)),$$
given by $ \chi \mapsto L(\chi)$, is bijective.\label{bijection1}
\item Let $L$ be a finite-dimensional simple left
$U(\D,\lambda)$-module, and $\chi \in \widehat{\G}$. Then $L \cong
L(\chi)$ if and only if there is a nonzero element $m \in L$ such
that $x_j \cdot m =0$, for all $ j \in I^+$, and $g \cdot m =
\chi(g)m$ for all $g \in \G$.\label{Labstract1}
\end{enumerate}
\end{Theorem}
\pf For all $1 \leq i \leq n$ we define
$$z'_i=z_{t(i)}, \eta'_i=\eta_{t(i)},u'_i=u_{t(i)}, g'_i=g_{t(n+i)},\chi'_i=\chi_{t(n+i)},a'_i=a_{t(n+i)}.$$
Let $W'=\oplus_{1 \leq i \leq n} ku'_i \in {\YDL},V'=\oplus_{1
\leq i \leq n} ka'_i \in {\YDG},$ and
$$U'=\mathfrak{B}(W') \# k[\Lambda],A'=\mathfrak{B}(V') \# k[\G].$$
We have checked in the proof of Theorem \ref{twistingbiproduct}
that
\begin{equation*}
\varphi (z'_i) = (\chi'_i)^{-1} \text{ and }\eta'_i(z) =
\varphi(z)(g'_i) \text{ for all } 1 \leq i \leq n, z \in \Lambda.
\end{equation*}
Thus we are in the non-degenerate situation of Section
\ref{computations} described in \eqref{conditionvarphi}. Let
$\sigma'$ be the corresponding 2-cocycle. Then by Theorem
\ref{maingeneric} the map
$$\{(\rho,\chi) \in \widehat{\Lambda} \times \widehat{\G} \mid (\rho,\chi) \text{ dominant}\} \to \Irr((U' \otimes A')^{\sigma'}),$$
given by $(\rho,\chi) \mapsto [L_{H'}(\rho,\chi)]$, is bijective,
where $H'=(U' \# A')^{\sigma'}$.

Note that $\eta'_i(z'_i)=\eta_{t(i)}(z_{t(i)})=
\chi_{t(i)}(g_{t(i)})=q_{t(i),t(i)}.$ Hence a pair $(\rho,\chi)$
is dominant if
\begin{equation}\label{meaningdom}
(\eta'_i(z'_i))^{m_i}\rho(z'_i)\chi(g'_i)=q_{t(i),t(i)}^{m_i}
\rho(z_{t(i)})\chi(g_{t(n+i)})=1
\end{equation}
for all $1 \leq i \leq n$, where the $m_i$ are integers $\geq 0$.

Then it follows from Theorem \ref{nondegenerate} that also
$$\{(\rho,\chi) \in \widehat{\Lambda} \times \widehat{\G} \mid (\rho,\chi) \text{ dominant}\} \to \Irr((U\otimes A)^{\sigma}),$$
given by $(\rho,\chi) \mapsto [L(\rho,\chi)]$, is bijective.

From \cite[Lemma 4.3]{Pairs} and the discussion above we finally obtain that the map
$$\{\chi \in \widehat{\G} \mid \chi \text{ dominant }\} \to \Irr(U(D,\lambda)),$$
given by $ \lambda \mapsto L(\lambda)$, is bijective, where $\chi
\in \widehat{\G}$ dominant means that the pair $(\rho,\chi)$ with
$\rho(z_i)=\chi(g_i), i \in I^-$, is dominant. By
\eqref{meaningdom} this latter condition says that
$q_{t(i),t(i)}^{m_i} \chi(g_{t(i)}g_{t(n+i)})=1$, that is, $\chi$
is dominant in the sense of our definition. This proves
\eqref{bijection1}, and \eqref{Labstract1} follows from Theorem
\ref{modules} \eqref{Labstract}. \epf

For all $1 \leq i \leq 2n$ let $g'_i=g_{t(i)},\chi'_i=\chi_{t(i)},
a'_{ij}=a_{t(i),t(j)}$. Then
$$\D'=\D((g'_i)_{1 \leq i \leq 2n}, (\chi'_i)_{1 \leq i \leq 2n},(a'_{ij})_{1 \leq i \leq 2n})$$ is a Cartan datum of finite type. Note that for all $1 \leq i,j \leq n$
\begin{equation}\label{a'}
a'_{ij}=a'_{n+i,n+j}, \text{ and } a'_{i,n+j}=0=a'_{n+i,j} \text{
for all } 1 \leq i,j \leq n,
\end{equation}
by Lemma \ref{linkingLemma} \eqref{Cartanlinkable}, and since
$t(i) \nsim t(n+j),t(j) \nsim t(n+i)$. Let $I'=\{1,2,\dots,2n\}$.
If $i,j \in I'$, then $t(i) \nsim t(j)$ in $I$ implies that $i
\nsim j $ in $I'$, but in general the converse is not true. We
define a family of linking parameters $\lambda'$ for $\D'$ by
$$\lambda'_{ij}=\begin{cases} \lambda_{t(i),t(j)} & \text{if $t(i) \nsim t(j)$ in $I$},\\
0 & \text{otherwise}\end{cases},$$ for all $1 \leq i,j \leq 2n,i
\nsim j$ in $I'$.

There is a surjective Hopf algebra map
\begin{equation*}
\pi : U(\D,\lambda) \to U(\D',\lambda'),
\end{equation*}
given by
$$\pi(x_k)=\begin{cases}
x_{i},&\text{ if } k=t(i),1 \leq i \leq 2n,\\
0,& \text{ otherwise}
\end{cases}$$
and $\pi(g)=g$ for all $g \in \G$. Note that $\pi$ preserves the
linking relations of $U(\D,\lambda)$, since $\lambda_{kl}=0$ for
all $k,l$ with $k \not\in t(I')$ or $l \not\in t(I')$; and $\pi$
preserves the Serre relations by \eqref{a'}, and since for all
$i,j \in I'$ with $t(i) \sim t(j)$ in $I$  and $i \nsim j$ in $I'$
the relation $\ad_c(x_i)(x_j) =0$ follows from the linking
relations in $U(\D',\lambda')$.

\begin{Cor}\label{redquotient}
The Hopf algebra map $\pi: U(\D,\lambda) \to U(\D',\lambda')$
induces a bijection
$$\pi^* : \Irr(U(\D',\lambda')) \to \Irr(U(\D,\lambda)).$$
\end{Cor}
\pf The isomorphism of Theorem \ref{twistingbiproduct} defines an
isomorphism
$$\Psi' : U(\D',\lambda') \xrightarrow{\cong} (U' \otimes A')^{\sigma}/(z_i \otimes g_i^{-1} - 1 \otimes 1 \mid i \in I^-),$$
and the claim follows from the proof of Theorem \ref{mainD}. \epf

\begin{Ex}\label{Excircle}
Suppose that $\D$ is a datum of finite Cartan type such that the
Dynkin diagram of the Cartan matrix of $\D$ is the disjoint union
of an even number of components of type $A_{n_l}, n_l \geq 2.$
Suppose that $\lambda$ is a family of linking parameters such that
the connected components are linked in a circle, the end of one
$A_{n_l}$ being linked to the beginning of the next $A_{n_{l+1}}.$
Then the Dynkin diagram of $\D'$ in Corollary \ref{redquotient} is
a union of components of type $A_1,$ and by Theorem \ref{mainD}
and Corollary \ref{quantumlinear2} we have an explicit description
of the finite-dimensional simple $U(\D,\lambda)$-modules.
\epf\end{Ex}

\subsection{Reduced data of Cartan type}\label{sectionreduced}

By Corollary \ref{redquotient} the computation of the
finite-dimensional simple  $U(\D,\lambda)$-modules can be reduced
to the non-degenerate case of the $U(\D',\lambda')$-modules. To
describe this case we introduce reduced data.

\begin{Def}
Let $\G$ be an abelian group, $n\geq 1,L_i,K_i\in \G, \chi_i\in
\widehat{\G}$ for all $1 \leq i \leq n,$ and $(a_{ij})_{1 \leq i,j
\leq n}$ a Cartan matrix of finite type. We say that
$\D_{red}=\D_{red}(\G, (L_i)_{1 \leq i \leq n}, (K_i)_{1 \leq i
\leq n}, (\chi_i)_{1 \leq i \leq n}, (a_{ij})_{1 \leq i,j \leq
n})$ is a {\em reduced datum of finite Cartan type} if for all
$1\leq i,j \leq n$
\begin{align}
\chi_i(K_j) \chi_j(K_i) &= \chi_i(K_i)^{a_{ij}}, \label{Cartanreduced1}\\
\chi_i(L_j) &= \chi_j(K_i),\label{Cartanreduced2}\\
K_iL_i &\neq 1, \chi_i(K_i) \neq 1.
\end{align}
A reduced datum $\D_{red}$ is called generic if $\chi_i(K_i)$ is
not a root of unity for all $1 \leq i \leq n$.
\end{Def}
\begin{Def}\label{Unondegenerate}
Let $\D_{red}$ be a reduced datum of finite Cartan type, and $X
\in {\YDG}$ with basis $x_1,\dots,x_n,y_1,\dots,y_n$, where $x_i \in
X_{L_i}^{\chi_i^{-1}}, y_{i}\in X_{K_i}^{\chi_i}$ for all $1 \leq
i \leq n.$ Let $l = (l_i)_{1 \leq i \leq n}$ be a family of
nonzero elements in $k$. Then we define $U(\D_{red},l)$ as the
quotient Hopf algebra of the smash product $k\langle
x_1,\dots,x_{n},y_1,\dots,y_n \rangle \# k[\G]$ modulo the ideal
generated by
\begin{align}
&\ad_c(x_i)^{1-a_{ij}}(x_j) \text{ for all } 1 \leq i,j \leq n, i \neq j,\label{reducedSerrex}\\
&\ad_c(y_i)^{1-a_{ij}}(y_j) \text{ for all } 1 \leq i,j \leq n, i \neq j,\label{reducedSerrey}\\
&x_iy_j - \chi_j(L_i) y_jx_i - \delta_{ij}l_i(1 - K_iL_i) \text{
for all } 1 \leq i,j \leq n.\label{reducedlinking}
\end{align}
For all $ 1 \leq i,j \leq n$ we let $q_{ij} = \chi_j(K_i)$.
\end{Def}

The discussion below shows that $U(\D',\lambda') \cong
U(\D_{red},l)$, where $\D'$ and $\lambda'$ are defined in the end
of the last section, and where
$$K_i =g'_i,L_i=g'_{n+i},\chi_i={\chi'}_{n+i},l_i =\lambda'_{i,i+n},a_{ij}=a'_{ij}$$
for all $1 \leq i,j \leq n$.

The definition of $U(\D_{red},l)$ is a special case of Definition
\ref{DefU}. Indeed define
$$g_i=L_i,g_{n+i}=K_i, \overline{\chi}_{i}=\chi_i^{-1}, \overline{\chi}_{n+i}=\chi_i \text{ for all } 1 \leq i \leq n,$$ and let $(a_{ij})_{1\leq i,j \leq 2n}$ be the diagonal block matrix consisting of two identical blocks $(a_{ij})_{1 \leq i,j \leq n}$ on the diagonal. Then
$$\D=\D(\G,(g_i)_{1 \leq i \leq 2n}, (\overline{\chi}_i)_{1 \leq i \leq 2n}, (a_{ij})_{1 \leq i,j \leq 2n})$$ is a datum of finite Cartan type. Note that \eqref{Cartanreduced1} and \eqref{Cartanreduced2} together are the Cartan condition \eqref{Cartan} for $\D$. We define a family $\lambda=(\lambda_{ij})_{1 \leq i,j \leq 2n,i \nsim j}$ of linking parameters for $\D$ for all $1 \leq i<j \leq 2n,i \not\sim j,$ by
$$\lambda_{ij}=\begin{cases} l_i & \text{if $1 \leq i \leq n,j=n+i$},\\
0 & \text{otherwise}\end{cases}.$$ The remaining values of
$\lambda$ are determined by \eqref{linkingnumbering}. Thus the
Dynkin diagram of $\D$ consists of two copies of the Dynkin
diagram of $(a_{ij})_{1 \leq i,j \leq n},$ and each vertex is
linked with its copy. Then
$$U(\D_{red},l) = U(\D,\lambda), \text{ where }y_i=x_{n+i} \text{ for all }1 \leq i \leq n,$$
since the set of relations \eqref{reducedSerrex},
\eqref{reducedSerrey}, \eqref{reducedlinking} coincides with the
set \eqref{Serrerelations}, \eqref{linkingrelations}.

The linked pairs of vertices of $\D$ are $(i,n+i)$ and $(n+i,i), 1
\leq i \leq n$. Hence we can apply Theorem \ref{twistingbiproduct}
to $\D$ with
$$I^-= \{1,2,\dots,n\},I^+=\{n+1,n+2,\dots,2n\}$$
and  $t=\id$. We call a character $\chi \in \widehat{\G}$ {\em
dominant} for $\D_{red}$ if there are natural numbers $m_i \geq 0,
1 \leq i \le n$, such that
\begin{equation}\label{dominantreduced}
\chi(K_iL_i) = q_{ii}^{m_i} \text{ for all } 1 \leq i \leq n.
\end{equation}

Then a character $\chi \in \widehat{\G}$ is dominant for
$\D_{red}$ if and only $\chi$ is dominant for $\D, \lambda$ and
$I^+$ since for all $1 \leq i \leq n$
$$\overline{\chi}_i(g_i) = \chi_i^{-1}(L_i)=\chi_i^{-1}(K_i) = q_{ii}^{-1}.$$

For $\chi \in \widehat{\G}$ we let $L(\chi)$ be the
left module over $U(\D_{red},l)=U(\D,\lambda)$ of Theorem
\ref{mainD}.

\begin{Cor}\label{mainreduced}
Let
$$\D_{red}=\D_{red}(\G, (L_i)_{1 \leq i \leq n}, (K_i)_{1 \leq i \leq n}, (\chi_i)_{1 \leq i \leq n}, (a_{ij})_{1 \leq i,j \leq n})$$
be a reduced, generic datum of finite Cartan type, and $l=(l_i)_{1
\leq i \leq n}$ be a family of nonzero elements in $k$.
\begin{enumerate}
\item The map
$$\{\chi \in \widehat{\G} \mid \chi \text{ dominant for }U(\D_{red},l)\} \to \Irr(U(D_{red},l)),$$
given by $ \chi \mapsto L(\chi)$, is bijective.
\item Let $L$ be a finite-dimensional simple left
$U(\D_{red},l)$-module, and $\chi \in \widehat{\G}$ a dominant
character. Then $L \cong L(\chi)$ if and only if there is a
nonzero element $m \in L$ such that $E_i \cdot m =0$, for all $ 1
\leq i \leq n$, and $g \cdot m = \chi(g)m$ for all $g \in
\G$.
 \item Let $\chi \in \widehat{\G}$ a
dominant character. Then
$$L(\chi)= \bigoplus_{\substack{1 \leq i_1,\dots,i_t \leq n,\\t \geq 0}} \ L(\chi)^{\chi_{i_1}^{-1} \cdots\chi_{i_t}^{-1}\chi}.$$
\end{enumerate}
\end{Cor}
\pf (1) and (2) follow from Theorem \ref{mainD}, and (3) from
Proposition \ref{actL} b). \epf

To see that $U_q(\mathfrak{g}), \mathfrak{g}$ a semisimple Lie
algebra, and other quantum groups in the literature are special
cases of $U(\D_{red}, l)$, we formulate the relations of
$U(\D_{red}, l)$ in terms of the usual generators $E_i$ and $F_i$.

\begin{Lem}\label{EF}
Let
$$E_i=y_i,F_i=x_iL_{i}^{-1},1 \leq i,j \leq n.$$
For each connected component $J \in \mathcal{X}$ we choose $q_J
\in k$ such that \eqref{DJ} holds. For all $1 \leq i,j \leq n$ let
\begin{equation*}
p_{ij}=\begin{cases} q_{ij}q_J^{-d_ia_{ij}},& \text{ if } i,j \in J, J \in \mathcal{X},\\
q_{ij},& \text{ if } i \nsim j.
\end{cases}
\end{equation*}
Then the relations \eqref{reducedSerrey}, \eqref{reducedSerrex}
and \eqref{reducedlinking} of $U(\D_{red},l)$ can be reformulated
in $U(\D_{red},l)$ for all $1 \leq i,j \leq n$ as
\begin{align}
&\sum_{s=0}^{1-a_{ij}}(-p_{ij})^s  \begin{bmatrix}
1-a_{ij}\\
s
\end{bmatrix}_{q_J^{d_i}} E_i^{1-a_{ij}-s} E_j E_i^s =0, i \neq j,i \in J, J \in \mathcal{X},\label{SerreE}\\
&\sum_{s=0}^{1-a_{ij}}(-p_{ij}^{-1})^s  \begin{bmatrix}
1-a_{ij}\\
s
\end{bmatrix}_{q_J^{d_i}} F_i^{1-a_{ij}-s} F_j F_i^s =0,  i \neq j,i \in J, J \in \mathcal{X},\label{SerreF}\\
&E_iF_j - F_jE_i = \delta_{ij} q_{ii}^{-1} l_i(K_i -
L_i^{-1}).\label{reducedlinking1}
\end{align}
The action of $\G$ is given for all $g \in \G$ and all $1 \leq i
\leq n$ by
\begin{align}
gE_ig^{-1} &= \chi_i(g) E_i, gF_ig^{-1} = \chi_i^{-1}(g) F_i,
\intertext{and the comultiplication by} \Delta(E_i) &= K_i \otimes
E_i + E_i \otimes 1, \Delta(F_i) = 1 \otimes F_i + F_i \otimes
L_i^{-1}.
\end{align}
\end{Lem}
\pf This follows from the defining relations of $U(\D_{red},l)$
and Lemma \ref{equivalentSerre}. \epf

\begin{Rem}\label{choose}
In the situation of Corollary \ref{mainreduced} let $\mathcal{M}$
be the set of all $n$-tuples $m=(m_i)_{1 \leq i \leq n}$ of
integers $\geq0$ for which there exists a dominant character
$\chi$ for $\D_{red}$ such that $\chi(K_iL_i)=q_{ii}^{m_i}, 1 \leq i \leq n$. For each $m=(m_i) \in \mathcal{M}$ we
choose a dominant character $\chi_m$ for $\D_{red}$ satisfying
$$\chi_m(K_iL_i)=q_{ii}^{m_i} \text{ for all }1 \leq i \leq n.$$
Then any dominant character $\chi$ for $\D_{red}$ has the form
\begin{equation}\label{uniqueform}
\chi = \psi \chi_m, m \in \mathcal{M},\psi \in \widehat{\G} \text{
with } \psi(K_iL_i)=1 \text{ for all } 1 \leq i \leq n,
\end{equation}
where $m$ and $\psi$ are uniquely determined. Note that the
algebra maps $U(\D_{red},l) \to k$ are all of the form
$\widetilde{\psi}$ where $\psi \in \widehat{\G}$ is a dominant
character satisfying
\begin{equation}\label{one}
\psi(K_iL_i)=1 \text{ for all } 1 \leq i \leq n,
\end{equation}
and where $\widetilde{\psi}(E_i)=0=\widetilde{\psi_i}(F_i)$ for
all $1 \leq i \leq n$, and $\widetilde{\psi}(g) =\psi(g)$ for all
$g \in \G$. Thus the one-dimensional $U(\D_{red},l)$-modules are
the modules $L(\psi)=k_{\widetilde{\psi}}$ where $\psi$ satisfies
\eqref{one}, and it follows from Corollary \ref{mainreduced} that
the finite-dimensional simple $U(\D_{red},l)$-modules are of the
form
\begin{equation}\label{generalcase}
L(\chi) \cong k_{\widetilde{\psi}} \otimes L(\chi_m),
\end{equation}
where $m \in \mathcal{M}$, and $\psi \in \widehat{\G}$  satisfies
\eqref{one}. \epf\end{Rem}

In the following remarks we discuss several special cases of
Corollary \ref{mainreduced}, namely the usual one-parameter
deformation $U_q(\mathfrak{g})$ of $U(\mathfrak{g})$, where
$\mathfrak{g}$ is a semisimple Lie algebra, Lusztig's version of
the one-parameter deformation with more general group-like
elements, and two-parameter deformations of $U(\mathfrak{gl_n})$
and $U(\mathfrak{sl_n})$. In each case the natural choice for the
characters $\chi_m$ is to define them by dominant weights of the
weight lattice of $\mathfrak{g}$.

\begin{Rem}\label{Jantzen}
The quantum group $U_q(\mathfrak{g})$ \cite[4.3]{J},
$\mathfrak{g}$  a semisimple Lie algebra, is a special case of
$U(\D_{red},l)$. Let $(a_{ij})_{1 \leq i,j \leq n}$ be the Cartan
matrix of $\mathfrak{g}$ with respect to a basis
$\alpha_1,\dots,\alpha_n$ of the root system of $\mathfrak{g}$,
$d_i \in \{1,2,3\}$ with $d_ia_{ij}=d_ja_{ji}$ for all $1 \leq i,j
\leq n$, and $0\neq q \in k$ such that $q^{2d_i} \neq 1$ for all $
1 \leq i \leq n$. Let $\G$ be a free abelian group with basis
$K_i,1 \leq i \leq n$, and $L_i=K_i$ or all $1 \leq i \leq n$.
Define the characters $\chi_i \in \widehat{\G}$ by
$$\chi_j(K_i) = q^{d_ia_{ij}} \text{ for all }1 \leq i,j \leq n.$$
Then $\D_{red}=\D_{red}(\G, (L_i), (K_i), (\chi_i), (a_{ij}))$ is
a reduced datum, and $$U(\D_{red},l) \cong U_q(\mathfrak{g}),$$
where $l_i=q^{2d_i}(q^{d_i}-q^{-d_i})^{-1}$ for all $1 \leq i \leq
n$. Note that in Lemma \ref{EF}, $p_{ij}=1$ for all all $1 \leq
i,j \leq n$.

Assume that $q$ is not a root of unity, and let $\chi \in
\widehat{\G}$. Let $\chi$  a dominant character for $\D_{red}$. By
definition there are integers $r_1,\dots,r_n \geq 0$ with
$\chi(K_i^2)=q^{2d_im_i}$, or equivalently
$$\chi(K_i)=\pm q^{d_im_i} \text{ for all }1 \leq i \leq n.$$
We say that $\chi$ is a character of type 1 if
$\chi(K_i)=q^{d_im_i}$ for all $1 \leq i \leq n$. Thus we see that
characters of type 1 correspond uniquely to dominant weights
$\lambda = \sum_{i=1}^n m_i \varpi_i$, where
$\varpi_1,\dots,\varpi_n$ are the fundamental weights
\cite[4.1]{J}. If the dominant weight $\lambda$ corresponds to the
character $\chi$, then
$$L(\lambda) \cong L(\chi),$$
where $L(\lambda)$ is the simple $U_q(\mathfrak{g})$-module
defined in \cite[5.5]{J}. This follows from Corollary
\ref{mainreduced} (2), since $L(\lambda)$ is a finite-dimensional
simple $U_q(\mathfrak{g})$-module and contains a non-zero element
$v_{\lambda}$ with
$$E_i \cdot v_{\lambda}=0,K_i \cdot v_{\lambda} = q^{(\lambda,\alpha_i)}v_{\lambda},1 \leq i \leq n,$$
where $(\lambda,\alpha_i)=d_im_i$ (see \cite[4.1]{J}) for all $
1\leq i \leq n$. \epf
\end{Rem}

\begin{Rem}\label{Lusztig}
Similarly the quantum group $\mathbf{U}$ defined in
\cite[3.1.1]{L}, where the Cartan datum $(I,\cdot)$ is of finite
type, is a special case of $U(\D_{red},l)$ (This also holds in
general, if we assume that $(a_{ij})$ is a symmetrizable
generalized Cartan matrix). In \cite[3.1.1]{L} a root datum
$(Y,X,<,>,...)$ of type $(I,\cdot)$ is given. Let $\G$ be the
group with generators $K_{\mu}, \mu \in Y$, and relations $K_0=1,
K_{\mu}K_{\mu'}=K_{\mu + \mu'}$  for all $\mu,\mu' \in Y$. Thus
$\G \cong Y$ is a free abelian group of finite rank. To avoid
confusion we denote the elements of $\G$ by $K'_{\mu}$ and not by
$K_{\mu}$. Define characters $\chi_i \in \widehat{\G}$ by
$$\chi_j(K'_{\mu})= v^{<\mu,j'>} \text{ for all } \mu \in Y, j \in
I.$$ For all $i,j \in I$ we define $d_i=\frac{i \cdot i}{2},
a_{ij}=<i,j'>$, and  $K_i=L_i=K'_{d_ii}$. Then $(a_{ij})$ is a
Cartan matrix of finite type, and
$$\chi_j(K_i) = v^{d_ia_{ij}} \text{ for all }i,j \in I.$$
Again $\D_{red}=\D_{red}(\G, (L_i), (K_i), (\chi_i), (a_{ij}))$ is
a reduced datum (where we identify $I$ with $\{1,2,\dots,n\}$),
and
$$U(\D_{red},l) \cong U_q(\mathfrak{g}),$$  where $l_i=v^{2d_i}(v^{d_i}-v^{-d_i})^{-1}$ for all $1 \leq i \leq n$.

Note that the matrix $(\chi_j(K_i))$ is symmetric as in the
previous remark,  but it is not assumed that the elements $K_i$
generate the group $\G$.

Let $M$ be a left $\mathbf{U}$-module. For any $ \lambda \in X$
let
$$M^{\lambda} = \{ m \in M \mid K'_{\mu}m = v^{<\mu,\lambda>}m \}$$
be the weight space of $M$ of weight $\lambda$. The category
$\mathcal{C}$ consists of all left $\mathbf{U}$-modules $M$ such
that
$$M=\oplus_{\lambda \in X} M^{\lambda}.$$
By \cite[3.5.5, 6.3.4]{L} the finite-dimensional simple
$\mathbf{U}$-modules in $\mathcal{C}$ are isomorphic to
$\Lambda_{\lambda}$, where $\lambda \in X$ is dominant, that is
$<i,\lambda> \in \mathbb{N}$ for all $i \in I$. Let $\lambda \in
X$ be dominant. The $\mathbf{U}$-module $\Lambda_{\lambda}$
contains a nonzero element $\eta_{\lambda}$ of weight $\lambda$
such that $E_i\cdot \eta_{\lambda}=0$ for all $i \in I$. By
Corollary \ref{mainreduced}, $\Lambda_{\lambda} \cong L(\chi)$,
where $\chi \in \widehat{\G}$ is defined by $\chi(K'_{\mu}) =
v^{<\mu,\lambda>}$ for all $\mu \in Y$. Note that $\chi(K_i) =
v^{d_im_i}$ for all $i \in I$, where $m_i=<i,\lambda> \in
\mathbb{N}$.

Thus  we recover the classification of finite-dimensional simple
$\mathbf{U}$-modules in the category $\mathcal{C}$. In Corollary
\ref{mainreduced} and \eqref{generalcase} this classification is
extended to {\em all} finite-dimensional simple
$\mathbf{U}$-modules. \epf\end{Rem}

\begin{Rem}\label{Benkart}
In \cite{BW} Benkart and Witherspoon give a classification of the
finite-dimensional simple modules over a two-parameter deformation
of the general linear and special linear Lie algebras. We will see
that their classification in \cite[Theorem 2.19]{BW} is a special
case of Corollary \ref{mainreduced}. Fix a natural number $n \geq
1$, nonzero elements $r,s \in k$, and assume that $rs^{-1}$ is not
a root of unity. Let $E$ be a euclidean vector space with inner
product $(,)$ and orthonormal basis
$\epsilon_1,\dots,\epsilon_{n+1}$. Define $\alpha_i= \epsilon_i -
\epsilon_{i+1}$ for all $1 \leq i \leq n$. Thus $\{\epsilon_i -
\epsilon_j \mid 1 \leq i,j \leq n+1, i \neq j \}$ is a root system
of type $A_n$ with basis $\alpha_1,\dots,\alpha_n$.

Let $a_1,\dots,a_{n+1},b_1,\dots,b_{n+1}$ be the basis of a free
abelian group of rank $2(n+1)$. Define characters $\chi_j \in
\widehat{\G}$ for all $1 \leq j \leq n$ by
$$\chi_j(a_i) = r^{(\epsilon_i,\alpha_j)},\chi_j(b_i)=s^{(\epsilon_i,\alpha_j)} \text{ for all } 1 \leq i \leq n+1.$$
Let $K_i=a_ib_{i+1},L_i=(a_{i+1}b_i)^{-1}$ and $q_{ij} =
\chi_j(K_i)$ for all $ 1 \leq i,j \leq n$. Then for all $1 \leq
i,j \leq n$,
\begin{align*}
\chi_j(K_i) &= r^{(\epsilon_i, \alpha_j)} s^{(\epsilon_{i+1},\alpha_j)}=\chi_i(L_j),\\ q_{ii}&=rs^{-1},\\
q_{ij}q_{ji} &= \begin{cases}
q_{ii}^{-1},& \text{ if } |i-j|=1,\\
1,& \text{ if } |i-j|>1.
\end{cases}
\end{align*}
Thus the Cartan condition \eqref{Cartan} is satisfied with the
Cartan matrix $(a_{ij})_{1 \leq i,j \leq n}$ of type $A_n$ defined
by
$$a_{ij} = \begin{cases}
-1,&  \text{ if } |i-j|=1,\\
0,& \text{ if } |i-j|>1,
\end{cases},$$
and $\D_{red}=\D_{red}(\G, (L_i), (K_i), (\chi_i), (a_{ij}))$ is a
reduced datum. Note that the matrix $(q_{ij})$ is not symmetric.
Choose $q \in k$ with $q^2=rs^{-1}$. Then the Serre relations
\eqref{SerreE} for $|i-j|=1$ are here
$$E_i^2E_{j}-p_{ij}(q+q^{-1})E_iE_{j}E_i+p_{ij}^2E_{j}E_i^2=0.$$
Since $q_{i,i+1}=s,q_{i+1,i}=r^{-1}$ we find
\begin{alignat*}{2}
&p_{i,i+1}(q+q^{-1})=r+s,&p_{i,i+1}^2&=rs,\\
&p_{i+1,i}(q+q^{-1})=r^{-1}+s^{-1},&\quad
p_{i+1,i}^2&=r^{-1}s^{-1},
\end{alignat*}
hence
\begin{align*}
&E_i^2E_{i+1}-(r+s)E_iE_{i+1}E_i+rsE_{i+1}E_i^2=0,\\
&E_iE_{i+1}^2-(r+s)E_{i+1}E_{i}E_{i+1}+rsE_{i+1}^2E_i=0.
\end{align*}
The Serre relations \eqref{SerreF} for $|i-j|=1$ are
\begin{align*}
&F_i^2F_{i+1}-(r^{-1}+s^{-1})F_iF_{i+1}F_i+r^{-1}s^{-1}F_{i+1}F_i^2=0,\\
&F_iF_{i+1}^2-(r^{-1}+s^{-1})F_{i+1}F_{i}F_{i+1}+r^{-1}s^{-1}F_{i+1}^2F_i=0.
\end{align*}
Thus we have established the relations (R6) and (R7) in \cite{BW},
and
$$U(\D_{red}, l) =U_{r,s}(\mathfrak{gl}_{n+1}),$$
defined in \cite{BW}, with $E_i=e_i,F_i=f_i$ for all $ 1 \leq i
\leq n$, where we take $l_i = rs^{-1}(r-s)^{-1}$ for all $1 \leq i
\leq n$, and $l=(l_i)$.

A dominant weight $\lambda$ is an element $\lambda$ in
$\mathbb{Z}\epsilon_1 \oplus \cdots \oplus \mathbb{Z}
\epsilon_{n+1}$ such that $(\alpha_i,\lambda) \in \mathbb{N}$ for
all $1 \leq i \leq n$. If $\lambda= \sum_{i=1}^{n+1}
\lambda_i\epsilon_i$, where $\lambda_i \in \mathbb{Z}$ for all $1
\leq i \leq n+1$, then $\lambda$ is dominant if and only if
$$\lambda_1 \geq \lambda_2 \geq \cdots \geq \lambda_{n+1}.$$
Following \cite{BW} we define for any dominant weight $\lambda$ a
character $\widehat{\lambda}$ by
$$\widehat{\lambda}(a_i) = r^{(\epsilon_i,\lambda)},\widehat{\lambda}(b_i) = s^{(\epsilon_i,\lambda)} \text{ for all } 1 \leq i \leq n+1.$$
Note that
$$\widehat{\lambda}(K_iL_i) =(rs^{-1})^{(\alpha_i,\lambda)}=q_{ii}^{(\alpha_i,\lambda)}\text{ for all } 1 \leq i \leq n+1.$$
Thus $\widehat{\lambda}$ is a dominant character in our sense. If
$m=(m_i)_{1 \leq i \leq n}$ is any family of natural numbers $\geq
0$, we define $\lambda_i=\sum_{k=i}^n m_k,1 \leq i \leq n$, and
$\lambda_{n+1}=0$. Then $m_i=(\alpha_i,\lambda)$ for all $1 \leq i
\leq n$, where $\lambda=\sum_{i=1}^{n+1} \lambda_i\epsilon_i$.

This shows that the classification of the finite-dimensional
simple $U_{r,s}(\mathfrak{gl}_{n+1})$-modules (and similarly
$U_{r,s}(\mathfrak{sl}_{n+1})$-modules) in \cite[Theorem 2.19]{BW}
is a special case of our Corollary \ref{mainreduced}.
\epf\end{Rem}

\subsection{The Finite-dimensional Case}\label{finite}

In this section we fix a finite abelian group $\G,$ a datum
$\mathcal{D} = \mathcal{D}(\Gamma, (g_i)_{1 \leq i \leq \theta},
(\chi_i)_{1 \leq i \leq \theta}, (a_{ij})_{1 \leq i,j \leq
\theta})$ of finite Cartan type, and a family $\lambda$ of linking
parameters for $\D.$

{\em We assume that for all $1 \leq i \leq \theta,$ $
\ord(q_{ii})$ is odd and $>3$, that $\ord(q_{ii})$ is prime to $3$ if
$i$ is in a component $G_2.$}

For each connected component $J \in \mathcal{X},$ and positive
root $\alpha$ of the root system $\Phi_J$ of $J$ let $x_{\alpha}$
be the root vector in the free algebra $k\langle x_j \mid j \in
J\rangle$ defined in \cite[Section 4.1]{AS2} generalizing the root
vectors in \cite{L}.
\begin{Def}\label{Defu}
Let $u(\D,\lambda)$ be the quotient Hopf algebra of  the smash
product $k\langle x_1,\dots,x_{\theta}\rangle \# k[\Gamma]$ modulo
the ideal generated by \eqref{Serrerelations},
\eqref{linkingrelations} and
\begin{align}\label{rootvectorrelations}
&x_{\alpha}^{N_J} \text { for all } \alpha \in \Phi_J^+, J \in
\mathcal{X}.&
\end{align}
\end{Def}

In addition we assume that {\em the linking graph of $\D,\lambda$
is bipartite}. By Lemma \ref{linkinggraph} this holds in
particular if $(a_{ij})$ is simply laced.

We proceed exactly as in the previous section and use the above
notations. The only difference is the definition of $\Lambda$  as
in \cite[proof of Theorem 5.17]{AS2}.

Let $\Lambda$ be the abelian group with generators $z_i,i \in
I^-,$ and relations $z_i^{n_i}=1$ for all $i \in I^-,$ where $n_i$
is the least common multiple of $\ord(g_i)$ and $\ord(\chi_i).$
For all $j \in I^-$ let $\eta_j \in \widehat{\Lambda}$ be defined
by $\eta_j(z_i) = \chi_j(g_i)$ for all $i \in I^-.$

Let $W \in {\YDL}$ with basis $u_i \in W_{z_i}^{\eta_i}, i \in
I^-,$ and $V \in {\YDG}$ with basis $a_j \in V_{g_j}^{\chi_j}, j
\in I^+.$ We define
$$U=\mathfrak{B}(W) \# k[\lambda]\text{ and }A=\mathfrak{B}(V) \# k[\lambda].$$
Then $U = \mathfrak{B}(W)$ and $A = \mathfrak{B}(V)$ by
\cite[Theorem 4.5]{AS2}.

As in the infinite case using Theorem \ref{groupcase} we define for all $i \in I^-$ the
algebra map $\gamma_i : A \to k$, the
$(\varepsilon,\gamma_i)$-derivation $\delta_i :A \to k$ and the
Hopf algebra map $\Phi : U \to A^{0\,cop}$. Let $\sigma$ be the 2-cocycle corresponding to $\Phi.$ As in
the proof of \cite[Theorem 5.17]{AS2} we obtain

\begin{Theorem}\label{twistingbiproductfinite}
The group-like elements $z_i \otimes g_i^{-1}, i \in {X}_2,$ are
central in $(U \otimes A)^{\sigma},$ and there is an isomorphism
of Hopf algebras
$$u(\D,\lambda) \cong (U \otimes A)^{\sigma}/(z_i \otimes g_i^{-1} - 1 \otimes 1 \mid i \in I^-)$$
mapping $x_i$  with $i \in I^-, x_j$ with $j\in I^+,$ resp. $g \in
\Gamma$ onto the residue classes of $u_i \otimes 1,1 \otimes a_j $
resp. $1 \otimes g.$
\end{Theorem}

For any $ \chi \in \widehat{\G}$ define $\rho \in
\widehat{\Lambda}$ by $\rho(z_i)=\chi(g_i)$ for all $i \in I^-.$
As above we define $L(\chi)=L(\chi,\rho)$ as a left module over
$U(\D,\lambda).$ As in the infinite case we define the reduced
datum $\D'$, and $u(\D',\lambda')$ by adding the root vector
relations, and the Hopf algebra projection $F : u(\D,\lambda) \to
u(\D',\lambda')$.

\begin{Theorem}\label{finitesimple}
The function
$$\widehat{\G} \to \Irr(u(\D,\lambda)), \; \chi \mapsto [L(\chi)],$$
is bijective, and pullback along $F$ defines a bijection
$$F^* : \Irr(u(\D',\lambda')) \to \Irr(u(\D,\lambda)).$$
\end{Theorem}
\pf This is shown as in the proof of Theorem \ref{mainD} and
Corollary \ref{redquotient} using Theorem \ref{dim1} (2).
 \epf  \medskip

In particular each $x_k, k \neq t(i), 1 \leq i \leq 2n,$ is
contained in the kernel of $F$, hence $x_k$ lies in the Jacobson
radical of $u(\D,\lambda)$ by Theorem \ref{finitesimple}. We give
another proof of this fact without using the bipartiteness
assumption.

\begin{Theorem}\label{Jacobson}
Let $\D,\lambda$ be as in the beginning of this section but where
the linking graph of $\D,\lambda$ is not necessarily bipartite.
Let $1 \leq k \leq \theta$ be a vertex which is not linked to any
other vertex. Then $x_k$ is contained in the Jacobson radical of
$u(\D,\lambda).$
\end{Theorem}
\pf Let $J$ be the connected component of the Dynkin diagram
containing $k$. Let $u_J$ be the subalgebra of $u(\D,\lambda)$
generated by all $x_j, j \in J,$ and let $u'$ be the subalgebra of
$u(\D,\lambda)$ generated by all $g \in \Gamma$ and $x_l, l
\not\in J.$ Then using the PBW-basis in \cite[Theorem 3.3]{Class}
it follows that  $$u(\D,\lambda) = u'u_J.$$ Since $i$ is not
linked, $x_k$ skew-commutes with all the generators of $u'.$ Hence
$u(\D,\lambda)x_ku(\D,\lambda)\subseteq u(\D,\lambda)u_J^+.$ Since
the augmentation ideal $u_J^+$ of $u_J$ is nilpotent we see that
$x_k$ generates a nilpotent ideal in $u(\D,\lambda).$ \epf

\end{document}